# IMPROVED KERNEL ESTIMATION OF COPULAS: WEAK CONVERGENCE AND GOODNESS-OF-FIT TESTING[1]

BY MAREK OMELKA, IRÈNE GIJBELS AND NOËL VERAVERBEKE

*Charles University Prague, Katholieke Universiteit Leuven and Hasselt University*

We reconsider the existing kernel estimators for a copula function, as proposed in Gijbels and Mielniczuk [*Comm. Statist. Theory Methods* **19** (1990) 445–464], Fermanian, Radulovič and Wegkamp [*Bernoulli* **10** (2004) 847–860] and Chen and Huang [*Canad. J. Statist.* **35** (2007) 265–282]. All of these estimators have as a drawback that they can suffer from a corner bias problem. A way to deal with this is to impose rather stringent conditions on the copula, outruling as such many classical families of copulas. In this paper, we propose improved estimators that take care of the typical corner bias problem. For Gijbels and Mielniczuk [*Comm. Statist. Theory Methods* **19** (1990) 445–464] and Chen and Huang [*Canad. J. Statist.* **35** (2007) 265–282], the improvement involves shrinking the bandwidth with an appropriate functional factor; for Fermanian, Radulovič and Wegkamp [*Bernoulli* **10** (2004) 847–860], this is done by using a transformation. The theoretical contribution of the paper is a weak convergence result for the three improved estimators under conditions that are met for most copula families. We also discuss the choice of bandwidth parameters, theoretically and practically, and illustrate the finite-sample behaviour of the estimators in a simulation study. The improved estimators are applied to goodness-of-fit testing for copulas.

**1. Introduction.** Consider a random vector $\mathbf{X} = (X_1, \ldots, X_d)^\mathsf{T}$ with joint cumulative distribution function $H$ and marginal distribution functions $F_1, \ldots,$

Received July 2008; revised October 2008.
[1]Supported by the IAP Research Network P6/03 of the Belgian State (Belgian Science Policy). This work was done while the first author was a postdoctoral researcher at the Katholieke Universiteit Leuven and the Universiteit Hasselt within the IAP Research Network. Support of the Research Project LC06024 is also highly appreciated.

$F_d$. According to Sklar's theorem [see, e.g., Nelsen (2006)], there exists a $d$-variate function $C$ such that

$$(1) \qquad H(x_1,\ldots,x_d) = C(F_1(x_1),\ldots,F_d(x_d)).$$

The function $C$ is called a copula, and it is, in itself, a joint cumulative distribution function on $[0,1]^d$ with uniform marginals. If the marginal distribution functions $F_1,\ldots,F_d$ are continuous, then the function $C$ is unique and $C(u_1,\ldots,u_d) = H(F_1^{-1}(u_1),\ldots,F_d^{-1}(u_d))$, where, for $j=1,\ldots,d$, $F_j^{-1}(u) = \inf\{x : F_j(x) \geq u\}$, with $u \in [0,1]$, is the quantile function of $F_j$. The copula $C$ "couples" the joint distribution function $H$ to its univariate marginals, capturing as such the dependence structure between the components of $\mathbf{X} = (X_1,\ldots,X_d)^\mathsf{T}$.

Methods for estimation of copulas usually depend on how much we are willing to assume about the joint distribution function $H$. In fully parametric approaches with parametric models for both the copula and the marginals, maximum likelihood estimation may be used. Nowadays semiparametric estimation is quite popular, in which one specifies a parametric copula and estimates the marginals nonparametrically. In this paper, we focus on nonparametric estimation of the copula making as such no restrictive distributional assumptions on the copula nor on the marginals.

For simplicity of the presentation, we will restrict to the case $d=2$, and consider an independent and identically distributed sample $(X_1,Y_1)^\mathsf{T},\ldots,(X_n,Y_n)^\mathsf{T}$ of a bivariate random vector $(X,Y)^\mathsf{T}$ with joint distribution function $H$ and marginal distribution functions $F$ and $G$.

Nonparametric estimation of copulas goes back to Deheuvels (1979) who proposed, in order to test for independence, the following empirical copula estimator:

$$C_n(u,v) = \frac{1}{n} \sum_{i=1}^n \mathbb{I}\{\hat{U}_i \leq u, \hat{V}_i \leq v\} \qquad \text{with } \hat{U}_i = F_n(X_i),\ \hat{V}_i = G_n(Y_i),$$

where $F_n$ and $G_n$ are the empirical cumulative distribution functions of the marginals, and where $\mathbb{I}\{A\}$ denotes the indicator of a set $A$. This estimator is asymptotically equivalent [up to a term $O(n^{-1})$] with the estimator based directly on Sklar's theorem given by

$$(2) \qquad C_n(u,v) = H_n(F_n^{-1}(u), G_n^{-1}(v))$$

with $H_n$ the empirical joint distribution function. Weak convergence studies of this estimator can be found in Gänssler and Stute (1987), Fermanian, Radulovič and Wegkamp (2004) and Tsukuhara (2005). Our Monte Carlo experiments showed that it is better to use the following (asymptotically equivalent) modification of the empirical copula:

$$C_n^{(\mathrm{E})}(u,v) = \frac{1}{n}\sum_{i=1}^{n}\mathbb{I}\{\hat{U}_i^{(\mathrm{E})} \leq u, \hat{V}_i^{(\mathrm{E})} \leq v\}$$

(3)
$$\text{with } \hat{U}_i^{(\mathrm{E})} = \frac{n}{n+1}F_n(X_i), \ \hat{V}_i^{(\mathrm{E})} = \frac{n}{n+1}G_n(Y_i),$$

which shifts the pseudo-observations $F_n(X_i)$ and $G_n(Y_i)$ a bit closer to the left corner of the unit interval $[0,1]$ [see, e.g., Genest, Ghoudi and Rivest (1995)].

Fermanian, Radulovič and Wegkamp (2004) also proposed a smoothed version of the empirical copula. Their proposal is a straightforward modification of (3) and the estimator is defined as

(4) $$\hat{C}_n^{(\mathrm{SE})}(u,v) = \hat{H}_n(\hat{F}_n^{-1}(u), \hat{G}_n^{-1}(v)),$$

where the quantities $\hat{H}_n$, $\hat{F}_n$ and $\hat{G}_n$ are given by

$$\hat{H}_n(x,y) = \frac{1}{n}\sum_{i=1}^{n} K_n(x - X_i, y - Y_i), \qquad \hat{F}_n(x) = \hat{H}_n(x, +\infty),$$

(5)
$$\hat{G}_n(x) = \hat{H}_n(+\infty, y)$$

with

$$K_n(x,y) = K\left(\frac{x}{b_n}, \frac{y}{b_n}\right), \qquad K(x,y) = \int_{-\infty}^{x}\int_{-\infty}^{y} k(s,t)\,ds\,dt,$$

where $k(s,t)$ is a given bivariate kernel density function, and $b_n$ is a bandwidth sequence tending to zero with $n$. Fermanian, Radulovič and Wegkamp (2004) proved weak convergence of this estimator.

There are two kernel type estimators in the literature that pay special attention to the correction of the boundary bias. This typical bias associated with kernel estimation is present since a copula has its support on the bounded set $[0,1]^2$. The first reference is the mirror-reflection type estimator originating from the work of Gijbels and Mielniczuk (1990) on copula density estimation. They take care of boundary bias correction through data-augmentation obtained by reflecting the original data with respect to the edges and the corners of the unit square. The second reference is the estimator of Chen and Huang (2007), who proposed to use a local linear kernel in order to deal with the bias near the boundaries of the unit square.

A first goal of the present paper is to prove the weak convergence of the estimators of Gijbels and Mielniczuk (1990) and Chen and Huang (2007) under the assumption that $C$ has bounded second order partial derivatives on $[0,1]^2$ (see Theorem 1 in Section 2). It turns out, however, that for many commonly-used families of copulas (e.g., Clayton, Gumbel, normal, Student), the latter condition is not satisfied and the bias behavior at the corners of the unit square precludes the weak convergence on the whole $[0,1]^2$. We therefore propose improved "shrinked" versions of the estimators of Gijbels and Mielniczuk (1990) and Chen and Huang (2007). This shrinking is done by including a weight function which removes the corner bias. In the same spirit, we also suggest a modification of the copula estimator (4) of Fermanian, Radulovič and Wegkamp (2004). In Theorem 2 we establish weak convergence for all newly proposed estimators. The finite-sample performance of the estimators is demonstrated via a simulation study. We discuss optimal bandwidth selection and compare the performances of the estimators using various well-known distance measures.

The second goal of the paper is to discuss the use of the various estimators of copulas in goodness-of-fit testing problems.

The paper is organized as follows. In Section 2, we introduce the improved kernel estimators and state the main theoretical results on weak convergence. In Section 3, we investigate the finite-sample performance of the newly-proposed estimators and compare these with performances of existing estimators. In Section 4, simulation results are reported for goodness-of-fit testing. The proofs of the weak convergence results are given in the Appendix.

**2. Nonparametric kernel estimators of a copula.**   In this section, we briefly discuss existing kernel estimators and propose important modifications. We also state the weak convergence results.

2.1. *Local linear kernel estimator.*   Chen and Huang (2007) constructed their estimator in the following way. In the first stage, they estimate marginals by

$$(6) \qquad \hat{F}_n(x) = \frac{1}{n}\sum_{i=1}^n K\left(\frac{x-X_i}{b_{n1}}\right), \qquad \hat{G}_n(y) = \frac{1}{n}\sum_{i=1}^n K\left(\frac{y-Y_i}{b_{n2}}\right)$$

with $K$ the integral of a symmetric bounded kernel function $k$ supported on $[-1,1]$. In the second stage, the pseudo-observations $\hat{U}_i = \hat{F}_n(X_i)$ and $\hat{V}_i = \hat{G}_n(Y_i)$ are used to estimate the joint distribution function of the unobserved $F(X_i)$ and $G(Y_i)$, which gives the estimate of the unknown copula $C$. To prevent boundary bias, Chen and Huang (2007) suggested using a local linear version of the kernel $k$ given by

$$(7) \qquad k_{u,h}(x) = \frac{k(x)\{a_2(u,h) - a_1(u,h)x\}}{a_0(u,h)a_2(u,h) - a_1^2(u,h)}\mathbb{I}\left\{\frac{u-1}{h} < x < \frac{u}{h}\right\},$$

where

$$a_l(u,h) = \int_{(u-1)/h}^{u/h} t^l k(t)\,dt \qquad \text{for } l = 0, 1, 2.$$

Finally, the local linear type estimator of the copula is given by

$$\hat{C}_n^{(\text{LL})}(u,v) = \frac{1}{n}\sum_{i=1}^n K_{u,h_n}\left(\frac{u-\hat{U}_i}{h_n}\right) K_{v,h_n}\left(\frac{v-\hat{V}_i}{h_n}\right), \tag{8}$$

where $K_{u,h}(x) = \int_{-\infty}^x k_{u,h}(s)\,ds$. Chen and Huang (2007) derived expressions for asymptotic bias, variance and mean squared error for this estimator and showed that a proper choice of the second stage smoothing constants $h = h_n$ may considerably decrease variance and mean squared error of the copula estimate. Moreover, their Monte Carlo experiments showed that the estimator $\hat{C}_n^{(\text{LL})}$ is quite insensitive to the choice of the constants $b_{1n}$ and $b_{2n}$ used for smoothing the marginals in the first stage. Variance considerations provided by the authors even showed that it is reasonable to take $b_{1n}$ and $b_{2n}$ as small as possible. Note that strong undersmoothing in the first stage, recommended in Chen and Huang (2007), results in using the pseudo-observations $(\hat{U}_i, \hat{V}_i)^\mathsf{T} = (\frac{2nF_n(X_i)-1}{2n}, \frac{2nG_n(Y_i)-1}{2n})^\mathsf{T}$, which is asymptotically equivalent to the mostly-used pseudo-observations defined in (3).

As already mentioned in the Introduction, the theoretical inconvenience of the estimator (8) is that for many common families of copulas (e.g., Clayton, Gumbel, normal, Student) the bias of the estimator at some of the corners of the unit square is only of order $O(h_n)$. As the optimal bandwidth for distribution function estimation is of order $O(n^{-1/3})$, this violates the $n^{1/2}$-order weak convergence on the whole $[0,1]^2$.

The problem is caused by unboundedness of second order partial derivatives of many copula families. Although parametric models with unbounded densities are rather rare in "standard" parametric models, copula families with unbounded densities are quite common. As a benchmark, we can take the normal bivariate density, which is usually supposed to be a well-behaved model. But the resulting normal copula density is unbounded.

To overcome this difficulty, we propose a method of shrinking the bandwidth when coming close to the borders of the unit square. The proposed method is based on the observation that, when calculating the bias of the estimator (8), we have to deal with terms of the form $h^2 C_{uu}(u,v)$, $h^2 C_{uv}(u,v)$ and $h^2 C_{vv}(u,v)$, where $C_{uu}(u,v)$, $C_{uv}(u,v)$ and $C_{vv}(u,v)$ are the second order partial derivatives of $C$; that is, $C_{uu}(u,v) = \partial^2 C(u,v)/\partial u^2$ and, similarly, for $C_{uv}(u,v)$, $C_{vv}(u,v)$. The problem is that, for many common families of copulas, these second order partial derivatives are not bounded, and, in fact,

# 6
M. OMELKA, I. GIJBELS AND N. VERAVERBEKE

a closer inspection of them shows that

(9)
$$C_{uu}(u,v) = O\left(\frac{1}{u(1-u)}\right), \quad C_{vv}(u,v) = O\left(\frac{1}{v(1-v)}\right),$$
$$C_{uv}(u,v) = O\left(\frac{1}{\sqrt{uv(1-u)(1-v)}}\right).$$

This is shown in Appendix D for Clayton, Gumbel, normal and Student copulas. In order to keep the bias bounded, we suggest an improved "shrinked" version of (8), which is given by

(10)
$$\hat{C}_n^{(\text{LLS})}(u,v) = \frac{1}{n}\sum_{i=1}^n K_{u,h_n}\left(\frac{u-\hat{U}_i}{b(u)h_n}\right)K_{v,h_n}\left(\frac{v-\hat{V}_i}{b(v)h_n}\right)$$
$$\text{with } b(w) = \min(\sqrt{w}, \sqrt{1-w}),$$

where the constant bandwidth $h_n$ is replaced by a bandwidth function $b(u)h_n$ that "shrinks" the value of the bandwidth close to zero at the corners of the unit square. A straightforward adaptation of the result of Chen and Huang (2007) gives that, for $(u/b(u), v/b(v)) \in [h_n, 1-h_n]^2$ (and no smoothing of the marginals in the first stage),

(11)
$$\text{Bias}\{\hat{C}_n^{(\text{LLS})}(u,v)\}$$
$$= \frac{\sigma_K^2}{2}h_n^2\{b^2(u)C_{uu}(u,v) + b^2(v)C_{vv}(u,v)\} + o(h_n^2),$$

(12)
$$\text{Var}\{\hat{C}_n^{(\text{LLS})}(u,v)\}$$
$$= \frac{1}{n}\text{Var}[\mathbb{I}\{U \leq u, V \leq v\} - C_u(u,v)\mathbb{I}\{U \leq u\} - C_v(u,v)\mathbb{I}\{V \leq v\}]$$
$$- \frac{h_n b_K}{n}[b(u)C_u(u,v)(1-C_u(u,v))$$
$$+ b(v)C_v(u,v)(1-C_v(u,v))] + o\left(\frac{h_n}{n}\right)$$

with $\sigma_K^2 = \int_{-1}^1 t^2 k(t)\,dt$, $b_K = 2\int_{-1}^1 tk(t)K(t)\,dt$ and $b(\cdot)$ as defined in (10). Taking $b(w) = 1$ gives back the bias and variance expressions for $\hat{C}_n^{(\text{LL})}$ in Chen and Huang (2007) (in case of no smoothing at the first stage).

The improvements are obtained by shrinking the bandwidth through the function $b(\alpha, w) = \min\{w^\alpha, (1-w)^\alpha\}$. Different choices of $\alpha$ or different choices of shrinking factors are possible, but our extensive investigations showed that $b(w) = \min\{\sqrt{w}, \sqrt{1-w}\}$ is overall a very good choice. The choice of a possible optimal shrinking factor is an open question.

2.2. *Mirror-reflection kernel estimator.* Another version of a kernel estimator for the copula might be obtained by integration of the estimator of the density of the copula introduced and studied in Gijbels and Mielniczuk (1990). This estimator deals with the boundary problem by the technique known as mirror-reflection. If a multiplicative kernel $k(x,y) = k(x)k(y)$ is used, then the mirror-reflection estimate of the copula has a simple form

$$\hat{C}_n^{(\mathrm{MR})}(u,v) = \frac{1}{n}\sum_{i=1}^{n}\sum_{\ell=1}^{9}\left[K\left(\frac{u-\hat{U}_i^{(\ell)}}{h_n}\right) - K\left(\frac{-\hat{U}_i^{(\ell)}}{h_n}\right)\right] \\ \times \left[K\left(\frac{v-\hat{V}_i^{(\ell)}}{h_n}\right) - K\left(\frac{-\hat{V}_i^{(\ell)}}{h_n}\right)\right], \quad (13)$$

where $\{(\hat{U}_i^{(\ell)}, \hat{V}_i^{(\ell)}), i=1,\ldots,n, \ell=1,\ldots,9\} = \{(\pm\hat{U}_i, \pm\hat{V}_i), (\pm\hat{U}_i, 2-\hat{V}_i), (2-\hat{U}_i, \pm\hat{V}_i), (2-\hat{U}_i, 2-\hat{V}_i), i=1,\ldots,n\}$.

The mirror-type estimator (13) faces the same "corner bias" problem as the local linear estimator (8). To prevent this problem, we can "shrink" the bandwidth similarly as in (10) and propose

$$\hat{C}_n^{(\mathrm{MRS})}(u,v) = \frac{1}{n}\sum_{i=1}^{n}\sum_{\ell=1}^{9}\left[K\left(\frac{u-\hat{U}_i^{(\ell)}}{b(u)h_n}\right) - K\left(\frac{-\hat{U}_i^{(\ell)}}{b(u)h_n}\right)\right] \\ \times \left[K\left(\frac{v-\hat{V}_i^{(\ell)}}{b(v)h_n}\right) - K\left(\frac{-\hat{V}_i^{(\ell)}}{b(v)h_n}\right)\right]. \quad (14)$$

2.3. *Transformation estimator.* The unboundedness of the densities of many copula families brings us back to Sklar's theorem in (1) and to the estimator (4) proposed in Fermanian, Radulovič and Wegkamp (2004).

To control the bias of this estimator in order to achieve weak convergence, we need the boundedness of the second order partial derivatives of the original joint distribution $H$. As the bivariate normal benchmark example shows, this condition may be considerably weaker than the requirement of the bounded second order derivatives of the underlying copula $C$.

A possible methodological objection to the estimator $\hat{C}_n^{(\mathrm{SE})}$, defined in (4), may be its dependence on the marginal distributions. This is confirmed by Monte Carlo simulations which show that, for a given copula, the success of this estimator depends on the marginals crucially.

As the copula function is invariant to increasing transformations of the margins, it is possible to transform the original data to $X_i' = T_1(X_i)$ and $Y_i' = T_2(Y_i)$, where $T_1$ and $T_2$ are increasing functions, and then use $(X_i', Y_i')$ instead of the original observations $(X_i, Y_i)$ in the estimator $\hat{C}_n^{(\mathrm{SE})}$. The aim of the transformation is to simplify the kernel estimation of the joint distribution. As the direct choice of functions $T_1$, $T_2$ is difficult, we propose the following procedure. Let us first construct the uniform pseudo-observations $\hat{U}_i^{(\mathrm{E})} = \frac{n}{n+1} F_n(X_i)$ and $\hat{V}_i^{(\mathrm{E})} = \frac{n}{n+1} G_n(Y_i)$. Then, for a given distribution function $\Phi$, put $\hat{S}_i = \Phi^{-1}(\hat{U}_i^{(\mathrm{E})})$ and $\hat{T}_i = \Phi^{-1}(\hat{V}_i^{(\mathrm{E})})$. Finally, use these transformed pseudo-observations $(\hat{S}_i, \hat{T}_i)$ instead of the original observations $(X_i, Y_i)$ in the estimator (5) of the joint distribution function. As we know, the marginals to be given by the function $\Phi$, the suggested estimator has, in the case of multiplicative kernel, the following simple formula:

$$
(15) \quad \hat{C}_n^{(\mathrm{T})}(u,v) = \frac{1}{n} \sum_{i=1}^n K\!\left(\frac{\Phi^{-1}(u) - \Phi^{-1}(\hat{U}_i^{(\mathrm{E})})}{h_n}\right) \\
\times K\!\left(\frac{\Phi^{-1}(v) - \Phi^{-1}(\hat{V}_i^{(\mathrm{E})})}{h_n}\right).
$$

The advantage of this estimator is that it is not affected by the marginal distributions. Further bias calculations show that, if we choose $\Phi$, such that $\frac{\Phi'(x)^2}{\Phi(x)}$ is bounded, we take care of the "corner bias problem" that is present if we try to estimate the joint distribution of pseudo-observations directly. The above condition is satisfied, for example, for $\Phi$ the normal cumulative distribution function.

2.4. *Main results.* The main theoretical contribution of this paper is the weak convergence of the kernel estimators $\hat{C}_n^{(\mathrm{LL})}$, $\hat{C}_n^{(\mathrm{LLS})}$, $\hat{C}_n^{(\mathrm{MR})}$, $\hat{C}_n^{(\mathrm{MRS})}$ and $\hat{C}_n^{(\mathrm{T})}$.

For notational convenience, let us denote $\hat{F}_n$ and $\hat{G}_n$ the estimates of the marginals that are used to construct pseudo-observations; that is, in the following we will write $\hat{U}_i = \hat{F}_n(X_i)$ and $\hat{V}_i = \hat{G}_n(Y_i)$. For the weak convergence results we need these functions to be asymptotically equivalent to the empirical cumulative distribution functions $F_n$, $G_n$; that is,

$$
(16) \quad \sup_x |\hat{F}_n(x) - F_n(x)| = o_p\!\left(\frac{1}{\sqrt{n}}\right), \\
\sup_y |\hat{G}_n(y) - G_n(y)| = o_p\!\left(\frac{1}{\sqrt{n}}\right),
$$

which further implies the standard weak convergence of the processes $\sqrt{n}(\hat{F}_n - F)$ and $\sqrt{n}(\hat{G}_n - G)$ to particular Brownian bridges. For technical reasons, we will also suppose that the functions $\hat{F}_n$ and $\hat{G}_n$ are nondecreasing, which excludes higher order kernels (taking negative values) for the estimation of the marginals.

It is easy to see that (16) is satisfied if we define pseudo-observations as $\hat{U}_i = \frac{2nF_n(X_i)-1}{2n}$, $\hat{V}_i = \frac{2nG_n(Y_i)-1}{2n}$, or in a way given in (3).

If we decide for kernel smoothing of the marginals given in (6), then it is well known [see, e.g., Lemma 7 of Fermanian, Radulovič and Wegkamp (2004)] that assumption (16) is met if there exists $\alpha > 0$ such that, uniformly in $x$,

$$F(x + b_{1n}) = F(x) + b_{1n}f(x) + o(b_{1n}^{1+\alpha}) \qquad \text{with } \sqrt{n}b_{1n}^{1+\alpha} \to 0,$$

where $f$ denotes the derivative of $F$ and, similarly, for $G$ involving $b_{2n}$.

Let $\mathbb{C}_n^{(\mathrm{LL})}$, $\mathbb{C}_n^{(\mathrm{LLS})}$, $\mathbb{C}_n^{(\mathrm{MR})}$, $\mathbb{C}_n^{(\mathrm{MRS})}$, $\mathbb{C}_n^{(\mathrm{T})}$ be suitably normalized empirical copula processes; that is, for $(u,v) \in [0,1]^2$,

$$\mathbb{C}_n^{(\cdot)} = \sqrt{n}\{C_n^{(\cdot)} - C(u,v)\}.$$

The proof of the following theorem is given in Appendix A. The termininology on stochastic processes (e.g., pinned C-Brownian sheet) is taken from Tsukuhara (2005). We refer the reader to this reference for details on the concepts used.

THEOREM 1. *Suppose that $H$ has continuous marginal distribution functions and that the underlying copula function $C$ has bounded second order partial derivatives on $[0,1]^2$. If $h_n = O(n^{-1/3})$ and (16) is satisfied, then the (kernel) copula processes $\mathbb{C}_n^{(\mathrm{LL})}$, $\mathbb{C}_n^{(\mathrm{MR})}$ converge weakly to the Gaussian process $G_C$ in $\ell^\infty([0,1]^2)$, having representation*

$$(17) \qquad G_C(u,v) = B_C(u,v) - C_u(u,v)B_C(u,1) - C_v(u,v)B_C(1,v),$$

*where $C_u$ and $C_v$ denote the first order partial derivatives of $C$, and $B_C$ is a two-dimensional pinned C-Brownian sheet on $[0,1]^2$; that is, it is a centered Gaussian process with covariance function*

$$(18) \qquad E[B_C(u,v)B_C(u',v')] = C(u \wedge u', v \wedge v') - C(u,v)C(u',v').$$

While Theorem 1 requires boundedness of the second order partial derivatives of the copula $C$, the weak convergence result of Fermanian, Radulovič and Wegkamp (2004) for the estimator $\mathbb{C}_n^{(\mathrm{SE})}$ given by (4) requires boundedness of the second order derivatives of the original joint distribution function $H$. This may or may not be more stringent, depending on the marginals. Unfortunately, Theorem 1 excludes many commonly-used families of copulas. The next theorem and Appendix D guarantee that the weak convergence of the proposed improved estimators $\mathbb{C}_n^{(\mathrm{LLS})}$, $\mathbb{C}_n^{(\mathrm{MRS})}$, $\mathbb{C}_n^{(\mathrm{T})}$ holds for commonly-used copulas such as Clayton, Gumbel, normal and Student copulas.

REMARK. A careful reader may find out that all the published weak convergence results for the empirical estimator (2) or the smoothed empirical estimator (4) require smoothness of the first order partial derivatives $C_u$ and $C_v$ of the copulas $C$ on $[0, 1]^2$. But this smoothness assumption usually is not true for the families which do not have bounded second order partial derivatives (e.g., Clayton, Gumbel, normal and Student). For instance, the first order partial derivatives of the Clayton copula are not continuous in the corner point $(0, 0)$. The second step of our proof given in Appendix B shows that it is sufficient to assume that

(19)    $C_u, C_v$ are continuous in $[0, 1]^2 \setminus \{(0, 0), (0, 1), (1, 0), (1, 1)\}$.

THEOREM 2.    *Suppose that $H$ has continuous marginal distribution functions and that the copula $C$ has bounded second order partial derivatives on $(0, 1)^2$ and satisfies (9) and (19). If $h_n = O(n^{-1/3})$ and (16) is satisfied, then the (kernel) copula processes $\mathbb{C}_n^{(\mathrm{LLS})}$ and $\mathbb{C}_n^{(\mathrm{MRS})}$ converge weakly to the Gaussian process $G_C$ in $\ell^\infty([0, 1]^2)$ given in Theorem 1.*

*Moreover, if the functions $\Phi'$ and $\frac{\Phi'(x)^2}{\Phi(x)}$ are bounded, then the above statement holds also for the process $\mathbb{C}_n^{(\mathrm{T})}$.*

## 3. Finite sample comparisons.

3.1. *Set up and performed comparisons.* In our simulation study, we always use the Epanechnikov kernel $k(x) = \frac{3}{4}(1 - x^2)\mathbb{I}\{|x| \leq 1\}$ and the bivariate multiplicative kernel $k(x, y) = k(x)k(y)$. The optimality of the Epanechnikov kernel in kernel density estimation was proven in Epanechnikov (1969). For background information on multivariate kernels see, for example, Wand and Jones (1995) and Fan and Gijbels (1996).

We investigate the performances of the estimators $C_n^{(\mathrm{E})}, \hat{C}_n^{(\mathrm{T})}, \hat{C}_n^{(\mathrm{LL})}, \hat{C}_n^{(\mathrm{MR})}$ and $\hat{C}_n^{(\mathrm{LLS})}$. We do not include the estimator $\hat{C}_n^{(\mathrm{SE})}$, defined in (4), because this estimator is too strongly influenced by the marginals, which makes the comparison difficult. For example, for a normal copula with normal marginals, the estimator $\hat{C}_n^{(\mathrm{SE})}$ usually does slightly better than its competitors. But, for a normal copula with, for example, exponential marginals, the performance of $\hat{C}_n^{(\mathrm{SE})}$ is considerably worse than its competitors. We do not present results for the modification of the mirror-type estimator $\hat{C}_n^{(\mathrm{MRS})}$ either, since its performance was found to be close to that of the estimator $\hat{C}_n^{(\mathrm{LLS})}$.

The performances of the various estimators were evaluated using two criteria: a Kolmogorov–Smirnov distance $\mathrm{KS}_n$ and a Cramér–von Mises distance

$CM_n$; that is,

$$KS_n = \sup_{u,v} |\hat{C}_n(u,v) - C(u,v)|,$$

$$CM_n = \sum_{i=1}^{n} [\hat{C}_n(\hat{U}_i, \hat{V}_i) - C(\hat{U}_i, \hat{V}_i)]^2,$$

where $\hat{C}_n$ stands for any of the investigated estimators, for example, $\hat{C}_n^{(E)}$. The corresponding statistics are denoted accordingly, for example, $KS_n^{(E)}$ and $CM_n^{(E)}$. Originally, we included the mean integrated asymptotic error $Q_n = n \iint [\hat{C}_n(u,v) - C(u,v)]^2 \, du \, dv$ as well. Not surprisingly, this measure behaves similarly to the Cramér–von Mises distance, since $CM_n \approx n \iint (\hat{C}_n - C)^2 \, dC$, but it is not so sensitive to the bias of the underlying copula estimator. See also Section 4.

For computational reasons, the supremum in the Kolmogorov–Smirnov distance $KS_n$ was replaced by a maximum over a grid of $101 \times 101$ points.

3.2. *Bandwidth choice*. The estimator $\hat{C}_n^{(LL)}$ involves bandwidths $b_{n1}$ and $b_{n2}$ (for estimation of the marginals) as well as a bandwidth $h_n$ when using local linear fitting to estimate the copula. Preliminary simulation results confirmed the results of Chen and Huang (2007), that the estimator $\hat{C}_n^{(LL)}$ [as well as its modification $\hat{C}_n^{(LLS)}$] cannot be improved by smoothing the marginals. Therefore, we simply work with the pseudo-observations $\hat{U}_i^{(E)} = \frac{n}{n+1} F_n(X_i)$ and $\hat{V}_i^{(E)} = \frac{n}{n+1} G_n(Y_i)$, with $F_n$ and $G_n$ the empirical cumulative distribution functions. This slightly differs from the strategy of strong undersmoothing recommended in Chen and Huang (2007), which more or less results in taking $\hat{U}_i = \frac{2nF_n(X_i)-1}{2n}$ and $\hat{V}_i = \frac{2nG_n(Y_i)-1}{2n}$. Nevertheless, the behavior of the resulting estimators is very similar.

For choosing the bandwidth $h_n$ for $\hat{C}_n^{(LL)}$ and $\hat{C}_n^{(LLS)}$, we rely on the expressions for asymptotic bias, variance and MISE derived in Chen and Huang (2007). From the main (Asymptotic) terms in (11) and (12), we derive the asymptotic mean squared error of the copula estimator in a given point $(u,v)$

(20) $\quad AMSE\{\hat{C}_n(u,v)\} = AVar\{\hat{C}_n(u,v)\} + [ABias\{\hat{C}_n(u,v)\}]^2.$

An optimal bandwidth is obtained by minimization of $\iint AMSE\{C_n(u,v)\} \, dC(u,v)$. As the true copula is unknown, this minimization cannot be carried out. A possible approach is then to consider a so-called reference copula. Chen and Huang (2007) proposed using a $t$-copula as a reference copula. But, as the second derivatives of the $t$-copula are not bounded, we experienced numerical difficulties and instabilities trying to apply this reference rule. We therefore decided to use Frank's copula, which has bounded second derivatives. The unknown parameter in Frank's copula family is estimated by inversion of Kendall's tau. The computational simplicity of this approach also makes the goodness-of-fit testing procedures, presented in Section 4, much more feasible.

Since the shrinkage of the bandwidth in the estimator $\hat{C}_n^{(\text{LLS})}$ removes the problem of possible unboundedness of the second order partial derivatives, there are plenty of families of copulas to use as a reference copula for this estimator. For simplicity and for more appropriate comparisons, we also use Frank's copula as a reference for $\hat{C}_n^{(\text{LLS})}$.

The asymptotic expansions (11) and (12) hold for the mirror-type kernel estimators $\hat{C}_n^{(\text{MR})}$ and $\hat{C}_n^{(\text{MRS})}$ as well; hence we also rely here on the same choice for $h_n$.

For the two improved estimators, a Frank copula based reference selection rule seems to give quite good performance (see Sections 3 and 4). A normal copula based reference rule tends to result in a too large bandwidth, whereas a Clayton copula based reference rule tends to give, on average, too small bandwidths. Thus Frank's reference rule seems to be a good compromise.

More problematic is the bandwidth choice for $\hat{C}_n^{(\text{T})}$, as we do not have asymptotic expressions for bias and variance here. We tried to minimize the expected mean squared integrated error

$$(21) \qquad \int_{-\infty}^{+\infty} \int_{-\infty}^{+\infty} [\hat{H}_n(x,y) - H(x,y)]^2 h(x,y)\,dx\,dy$$

taking a bivariate normal distribution $H$, with corresponding density $h$, as a reference distribution [see Jin and Shao (1999)]. However, the resulting bandwidth selector turned out to be too big. A possible explanation is that such a selection rule does not take into account that we rely on pseudo-observations $(\hat{U}_i, \hat{V}_i)$ instead of on the unobservable $(U_i, V_i)$. In our simulation study, we then used the above mentioned bandwidth divided by a factor two. This seems to be a reasonable ad-hoc solution. To further investigate the bandwidth selection problem, for $\hat{C}_n^{(\text{T})}$, we calculated the ratio of the bandwidth selected via (21) to the one selected via searching for a bandwidth that minimizes the criterion $\text{KS}_n(h)$ [resp., $\text{CM}_n(h)$] over a grid of $h$-values. Table 1 summarizes the obtained average ratios, for various values of Kendall's tau from 2000 simulated samples. Note that the ratios stay quite stable across different families of copulas as well as for different sample sizes. This suggests that it may be possible to find a reliable reference-based rule for $\hat{C}_n^{(\text{T})}$ as well.

The simulation studies reported below showed a promising performance for the transformation estimator $\hat{C}_n^{(\text{T})}$. A good bandwidth selection rule is missing, for the moment, and is subject of further research.

TABLE 1
*Average ratios of bandwidths for $\hat{C}_n^{(T)}$ selected from minimizing (21) and from criteria $\mathrm{KS}_n(h)$ and $\mathrm{CM}_n(h)$, respectively, for different Kendall's $\tau$ and sample sizes $n$*

|  | Clayton | | | | Frank | | | | Normal | | | |
| --- | --- | --- | --- | --- | --- | --- | --- | --- | --- | --- | --- | --- |
|  | $\tau=0.25$ | | $\tau=0.75$ | | $\tau=0.25$ | | $\tau=0.75$ | | $\tau=0.25$ | | $\tau=0.75$ | |
|  | $\mathrm{KS}_n$ | $\mathrm{CM}_n$ | $\mathrm{KS}_n$ | $\mathrm{CM}_n$ | $\mathrm{KS}_n$ | $\mathrm{CM}_n$ | $\mathrm{KS}_n$ | $\mathrm{CM}_n$ | $\mathrm{KS}_n$ | $\mathrm{CM}_n$ | $\mathrm{KS}_n$ | $\mathrm{CM}_n$ |
| $n=50$ | 1.20 | 1.28 | 1.53 | 2.00 | 1.23 | 1.38 | 1.53 | 2.04 | 1.28 | 1.37 | 1.26 | 1.62 |
| $n=150$ | 1.13 | 1.27 | 1.38 | 2.08 | 1.21 | 1.47 | 1.55 | 2.07 | 1.20 | 1.36 | 1.24 | 1.60 |

3.3. *Simulation results.* An extensive simulation study was carried out to compare the performances of all estimators using the performance measures $\mathrm{KS}_n$ (the Kolmogorov–Smirnov distance) and $\mathrm{CM}_n$ (the Cramér–von Mises distance). To illustrate our main findings, we only report on results obtained for the following two simulation models:

*Model* 1. Frank copula with Kendall's $\tau = 0.25$;

*Model* 2. Clayton copula with Kendall's $\tau = 0.75$.

Models 1 and 2 represent very different copula functions. The copula in model 1 has bounded second order partial derivatives and presents a case of mild dependence, whereas the copula in model 2 has unbounded second order partial derivatives and shows a strong dependence between $X$ and $Y$. From each model, we simulated 10,000 samples of sample size $n = 150$.

Figure 1 shows the boxplots of the performance measures $\mathrm{KS}_n$ and $\mathrm{CM}_n$ for model 1 (top panels) and model 2 (bottom panels). Note that, for model 1, the estimators $\hat{C}_n^{(LL)}$, $\hat{C}_n^{(MR)}$, $\hat{C}_n^{(LLS)}$ and $\hat{C}_n^{(T)}$ perform very comparable for the Cramér–von Mises distance measure. For the Kolmogorov–Smirnov performance measure, the estimator $\hat{C}_n^{(T)}$ performs slightly but significantly worse than the three other estimators. The latter is likely caused by the usage of a too small bandwidth, as can be anticipated by looking at the different results for the two performance measures and Table 1. For model 2 (bottom panels), one clearly sees a better performance of the improved estimators $\hat{C}_n^{(LLS)}$ and $\hat{C}_n^{(T)}$, especially when looking at the performance measure $\mathrm{CM}_n$. This is as to be expected, since the copula in model 2 has unbounded second order partial derivatives and the measure $\mathrm{CM}_n \approx n \iint (\hat{C}_n - C)^2 \, dC$ is most affected by points with higher values of the copula density $c(u,v)$ (which usually correspond with points with higher values of the second order partial derivatives $C_{uu}$ and $C_{vv}$). In other words, the performance measure $\mathrm{CM}_n$ is more sensitive to the corner bias problem than the measure $\mathrm{KS}_n$. For model 1, there was no need for "shrinking" the bandwidth since the copula function has bounded second order partial derivatives. Nevertheless, the "shrinked-bandwidth" (improved) local linear estimator also performs very well for this model. The very promising performance of the transformation estimator $\hat{C}_n^{(T)}$ in view of the Cramér–von Mises distance $CM_n$ may be partially understood by the fact that, in view of Table 1, our ad-hoc rule of bandwidth choice for the estimator $\hat{C}_n^{(T)}$ is almost optimal in that situation.

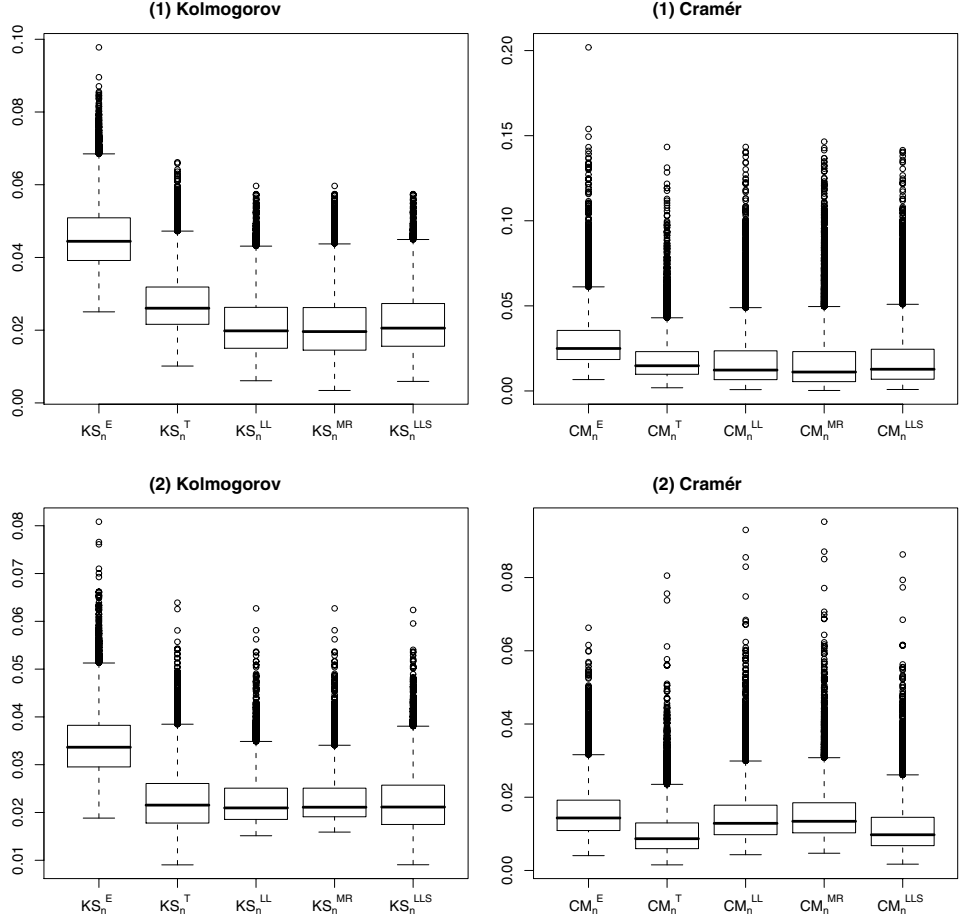

FIG. 1. *Boxplots of quantities* $KS_n$ *and* $CM_n$ *for different copula estimators. Top panels: model 1; Bottom panels: model 2.*

From the more extensive simulation study, we further report the following observations. All kernel copula estimators usually improve upon the empirical estimate $\hat{C}_n^{(E)}$. For copulas with bounded second order partial derivatives, the performances of the estimators $\hat{C}_n^{(LL)}$, $\hat{C}_n^{(MR)}$ and $C_n^{(E)}$ become very comparable, especially with increasing sample size. Overall, $\hat{C}_n^{(MR)}$ works slightly better for copulas with bounded second order partial derivatives [e.g., Frank, Farlie–Gumbel–Morgenstern, Ali–Mikhail–Haq; see Nelson (2006)] and mild dependence, with a significant improvement for copulas very close to independence copulas. On the other hand, the local linear kernel estimator $\hat{C}_n^{(\text{LL})}$ is preferable [compared to $\hat{C}_n^{(\text{MR})}$] in the remaining cases.

To gain further insights in the kernel estimators, we examined the dependence of these estimators on the bandwidth. We again use models 1 and 2 to illustrate our findings. For brevity, we present results only for the estimator $\hat{C}_n^{(\text{LL})}$, since similar findings can be reported on for the other kernel estimators.

Figure 2 illustrates the performance of the copula $\hat{C}_n^{(\text{LL})}$ with a fixed bandwidth $h$, in view of the performance measures $\text{KS}_n$ and $\text{CM}_n$, for models 1 and 2 (top and bottom panels, resp.). For comparison purposes, we also include (at the far left of the horizontal axis) the boxplot summarizing the results for the empirical copula $\hat{C}_n^{(\text{E})}$. In addition, we provide in each picture a (vertical) boxplot that indicates the bandwidths selected for $\hat{C}_n^{(\text{LL})}$ via (20).

Note that the effect of bandwidth choice is most noticeable from the Kolmogorov–Smirnov quantity $\text{KS}_n$. This is particularly true for the Clayton copula, model 2. For model 1, the estimator $\hat{C}_n^{(\text{LL})}$ improves upon the empirical copula $\hat{C}_n^{(\text{E})}$ for almost all $h$-values in the considered range of values. For model 2, however, which presents a case of stronger dependence, a kernel estimator comes with a gain, but only for a carefully selected bandwidth.

From the vertically displayed boxplots of bandwidths selected, we can further remark that a bandwidth selected via (20) works in fact quite satisfactory. This is particularly true in case of mild dependence and for copulas with bounded second derivatives (such as model 1). It may lead to a slight oversmoothing in a situation of strong dependence and for copulas with unbounded second derivatives (cf. model 2). It is worth mentioning though, that the presented results for model 2 are almost among the "worst-case" scenarios here.

**4. Goodness-of-fit tests for copulas.** When modelling multivariate data using copulas, a popular method is to estimate marginals nonparametrically and the copula in a parametric way. This requires choosing a suitable family of copulas for the data at hand, which is not an easy task. In this section, we focus on testing the null hypothesis

$$H_0 : C \in \mathcal{C}_0,$$

where $\mathcal{C}_0 = \{C_\theta, \theta \in \Theta\}$ is a given parametric family of copulas.

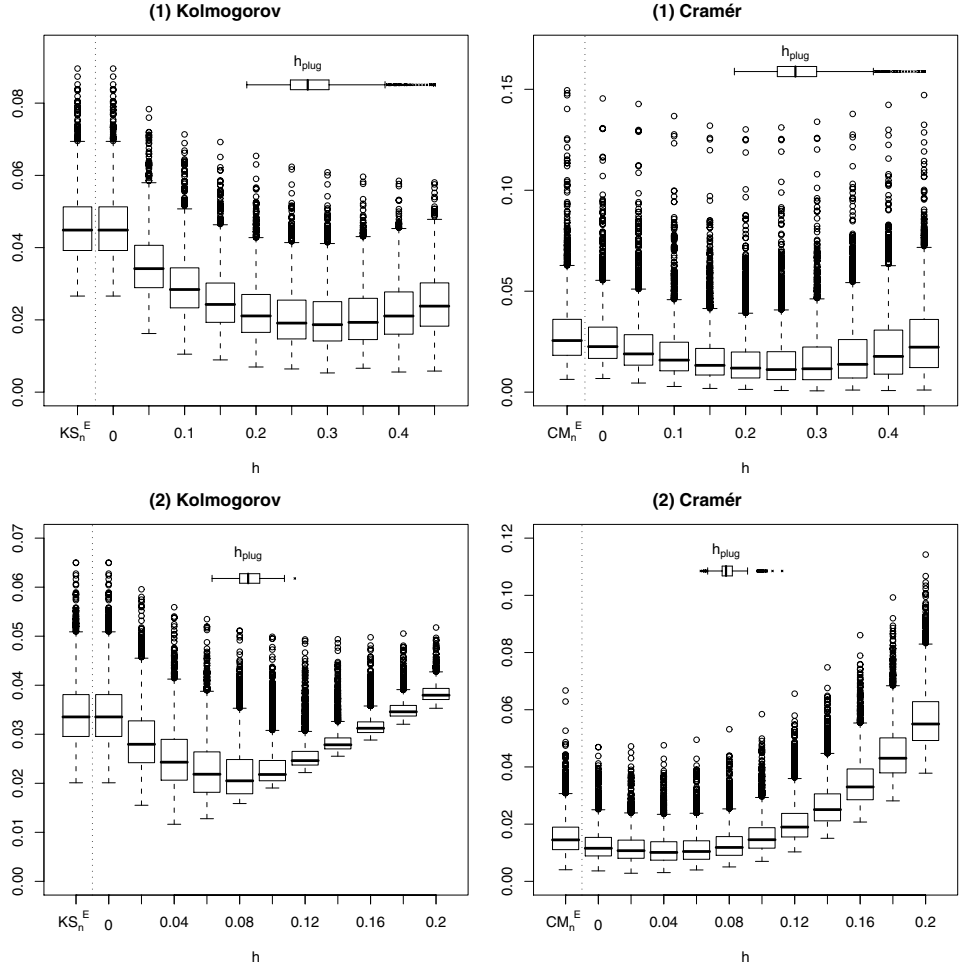

Fig. 2. *Boxplots of the quantities* $\mathrm{KS}_n$ *and* $\mathrm{CM}_n$ *for different values of fixed bandwidths for the estimator* $\hat{C}_n^{(\mathrm{LL})}$, *and boxplot (far left) for the estimator* $\hat{C}_n^{(\mathrm{E})}$. *Top panels: model 1; Bottom panels: model 2.*

Many testing methods have been proposed. See, for example, Chen and Huang (2007) and the review paper of Genest, Rémillard and Beaudoin (2008). The latter paper included a simulation study on classical goodness-of-fit measures such as the Kolmogorov–Smirnov and the Cramér–von Mises statistics, which we denote by (allowing a small abuse of previous notation)

$$\mathrm{KS}_n^{(\mathrm{E})} = \sup_{u,v} |C_n^{(\mathrm{E})}(u,v) - C_{\hat{\theta}_n}(u,v)|,$$

(22)

$$\mathrm{CM}_n^{(\mathrm{E})} = \sum_{i=1}^{n} [C_n^{(\mathrm{E})}(\hat{U}_i, \hat{V}_i) - C_{\hat{\theta}_n}(\hat{U}_i, \hat{V}_i)]^2,$$

where $\hat{\theta}_n$ is an estimate of the unknown parameter $\theta_0$ based on the inversion of the observed Kendall's tau.

The aim of this section is to investigate the size and power properties for testing procedures based on the test statistics $\mathrm{KS}_n^{(\mathrm{LL})}$, $\mathrm{KS}_n^{(\mathrm{LLS})}$, $\mathrm{CM}_n^{(\mathrm{LL})}$ $\mathrm{CM}_n^{(\mathrm{LLS})}$ computed by replacing $C_n^{(\mathrm{E})}$ in (22) with $\hat{C}_n^{(\mathrm{LL})}$ or $\hat{C}_n^{(\mathrm{LLS})}$. In addition, we consider here the test statistic

$$Q_n^{(\mathrm{E})} = \iint [C_n^{(\mathrm{E})}(u,v) - C_{\hat{\theta}_n}(u,v)]^2 \, du \, dv$$

and its $\hat{C}_n^{(\mathrm{LL})}$ and $\hat{C}_n^{(\mathrm{LLS})}$ versions. The double integral in the definition of $Q_n^{(\cdot)}$ was approximated by a double sum over a grid of $101 \times 101$ points.

Since the asymptotic distributions of these test statistics are too complex, a parametric bootstrap is used. This procedure runs as follows:

(1) By inversion of the empirical Kendall's $\tau$, estimate the unknown parameter $\theta$ of the null hypothesis family by $\hat{\theta}_n$ and compute the test statistic $\mathrm{KS}_n^{(\cdot)}$ [where the superscript $(\cdot)$ refers to any of the considered estimators of the copula];
(2) Generate $\{(U_i^*, V_i^*)\}_{i=1}^n$ from the copula $C_{\hat{\theta}_n}$ and use them as original observations to compute $\hat{\theta}_n^*$ and $\mathrm{KS}_n^{*(\cdot)}$;
(3) Repeat step (2) $B$-times;
(4) Estimate the $p$-values as

$$p_{\mathrm{KS}_n^{(\cdot)}} = \frac{1 + \#\{\mathrm{KS}_n^{*(\cdot)} \geq \mathrm{KS}_n^{(\cdot)}\}}{B+1}.$$

See Davison and Hinkley (1997).

For any of the other test statistics, we proceed similarly, replacing $\mathrm{KS}_n^{(\cdot)}$ by $\mathrm{CM}_n^{(\cdot)}$ or $Q_n^{(\cdot)}$.

According to Genest and Rémillard (2008), the validity of this bootstrap procedure requires the weak convergence of the copula processes $\mathbb{C}_n^{(\mathrm{LL})}$ and $\mathbb{C}_n^{(\mathrm{LLS})}$. For the latter process, the weak convergence is justified for all copula families considered in our simulation study, by Theorem 2. In contrast, only for Frank's copula the condition of Theorem 1 is satisfied when dealing with the weak convergence of the process $\mathbb{C}_n^{(\mathrm{LL})}$. Lemma C.1 of Appendix C shows that the test based on $\hat{C}_n^{(\mathrm{LL})}$ holds asymptotically the level even for the other families $\mathcal{C}_0$ appearing in the simulation study.

The setup of our simulation study closely follows that of Genest, Rémillard and Beaudoin (2008). The sample size is $n = 150$, and we take 999 number of bootstrap samples. Three values of Kendal's tau are considered, namely $\tau = 0.25, 0.50, 0.75$, for the following copula families: Clayton, Gumbel, Frank, normal and Student with four degrees of freedom (df). We use the R-computing environment, version 2.5.0 [see R Development Core Team (2007)], with copula package [see Yan (2007)]. For approximating the level of the test (i.e., under the null hypothesis) we use 6000 repetitions. The estimated powers of the test statistics are based on 1500 repetitions.

The results of the simulations are presented in Tables 2, 3 and 4. For ease of the reader, the estimated values for the size of the test statistics are presented in italics. Furthermore, for each testing problem, we highlighted the "best" power performances using bold characters. Readers should be aware of the fact that these estimated powers and sizes (using, resp., 1500 and 6000 repetitions) are of course subject to Monte Carlo approximation errors. A conservative upper bound (relying on a binomial distribution with parameters $B$ and $p$) for these approximations errors (in terms of standard deviation) is for the size estimates 0.28% (using $B = 6000$ and $p = 0.05$) and for the power estimates 1.29% (using $B = 1500$ and $p = 0.5$, for getting to an upper bound).

A summary of conclusions from the simulations results is as follows:

- The use of a kernel estimator [e.g., $\hat{C}_n^{(\text{LLS})}$] in goodness-of-fit testing seems to be promising in case true copulas are in the Clayton, Gumbel and Frank families, and we consistently improve upon the power for the Kolmogorov–Smirnov test;
- If a kernel estimator improves upon the power, it is most noticeable when the dependence is weaker, and is greatest for the Kolmogorov–Smirnov test;
- The power of the test statistics $Q_n$ is usually somewhere between the power of Kolmogorov–Smirnov and the Cramér–von Mises test statistics;
- The power of the test statistics based on the improved estimator $\hat{C}_n^{(\text{LLS})}$ is usually higher than this for test statistics based on $\hat{C}_n^{(\text{LL})}$ for alternatives with unbounded second order partial derivatives;
- For a true Frank copula and a Cramér–von Mises test statistic, the estimator $\hat{C}_n^{(\text{LLS})}$ is usually the best choice;
- The use of kernel estimators seems to be promising for Archimedean families of copulas (Clayton, Frank, Gumbel) but is somewhat questionable for elliptical families of copulas (normal, Student). Although kernel estimators may improve the power against Clayton and Gumbel alternatives, a loss in power is noticed for Frank alternatives. This holds in particular for $\hat{C}_n^{(\text{LLS})}$-based statistics.

## APPENDIX A: PROOF OF THEOREM 1

For simplicity, we will suppress the dependence on $n$ in the notation of pseudo-observations $(\hat{U}_i, \hat{V}_i)^\mathsf{T}$ and write, simply, $\hat{U}_i = \hat{F}_n(X_i) = \hat{F}_n(F^{-1}(U_i))$ and $\hat{V}_i = \hat{G}_n(Y_i) = \hat{G}_n(G^{-1}(V_i))$, where $(U_i, V_i)$ have a joint distribution function given by the copula $C$.

As our proof is a straightforward adaptation of the ideas used in van der Vaart and Wellner (2007), we would like to clarify one point. In the following, we will encounter the expectations of the form $\mathrm{E}g(\hat{U}, \hat{V})$, where $g$ is a measurable function on $[0,1]^2$ and $\hat{U} = \hat{F}_n(F^{-1}(U))$, $\hat{V} = \hat{G}_n(G^{-1}(V))$. In these types of expectations, the estimators of the marginal distribution functions $\hat{F}_n, \hat{G}_n$ are considered to be fixed (nonrandom) functions, and the expectation is taken only with respect to $(U, V)$ with joint distribution given by the copula $C$. Formally,

$$\mathrm{E}g(\hat{U}, \hat{V}) = \mathrm{E}_{U,V}[g(\hat{U}, \hat{V})|(X_1, Y_1), \ldots, (X_n, Y_n)],$$

whenever the integral on the right-hand side exists.

**A.1. Weak convergence of the process $\mathbb{C}_n^{(\mathrm{LL})}$.** In view of the previous remark, we decompose

$$
\begin{aligned}
&\sqrt{n}(\hat{C}_n^{(\mathrm{LL})}(u,v) - C(u,v)) \\
(23) \quad &= \frac{1}{\sqrt{n}}\left[\sum_{i=1}^n K_{u,h_n}\left(\frac{u-\hat{U}_i}{h_n}\right)K_{v,h_n}\left(\frac{v-\hat{V}_i}{h_n}\right) - C(u,v)\right] \\
&= A_n^{h_n}(u,v) + B_n(u,v) + C_n^{h_n}(u,v),
\end{aligned}
$$

where

$$
\begin{aligned}
A_n^{h_n}(u,v) = \frac{1}{\sqrt{n}}\bigg[&\sum_{i=1}^n K_{u,h_n}\left(\frac{u-\hat{U}_i}{h_n}\right)K_{v,h_n}\left(\frac{v-\hat{V}_i}{h_n}\right) \\
(24) \quad &- \mathbb{I}\{U_i \leq u, V_i \leq v\} \\
&- \mathrm{E}\left(K_{u,h_n}\left(\frac{u-\hat{U}}{h_n}\right)K_{v,h_n}\left(\frac{v-\hat{V}}{h_n}\right) - C(u,v)\right)\bigg]
\end{aligned}
$$

and

$$(25) \quad B_n(u,v) = \frac{1}{\sqrt{n}}\sum_{i=1}^n [\mathbb{I}\{U_i \leq u, V_i \leq v\} - C(u,v)],$$

$$(26) \quad C_n^{h_n}(u,v) = \sqrt{n}\mathrm{E}\left[K_{u,h_n}\left(\frac{u-\hat{U}}{h_n}\right)K_{v,h_n}\left(\frac{v-\hat{V}}{h_n}\right) - C(u,v)\right].$$

Our proof will be divided into two steps. First, we will show, in Step 1, that $\sup_{u,v}|A_n^{h_n}| = o_p(1)$. Then, we will prove, in Step 2, that

$$\sup_{u,v}|C_n^{h_n}(u,v) - \partial_1 C(u,v)\sqrt{n}[F_n^*(u) - u] - \partial_2 C(u,v)\sqrt{n}[G_n^*(v) - v]| = o_P(1),$$

TABLE 2
*Percentage of rejection of $H_0$ by various tests for samples of size $n = 150$ arising from different copula models with $\tau = 0.25$*

| Copula under $H_0$ | True copula | Cramér–von Mises | | | Kolmogorov–Smirnov | | | MISE | | |
|---|---|---|---|---|---|---|---|---|---|---|
| | | $\mathrm{CM}_n^{(\mathrm{E})}$ | $\mathrm{CM}_n^{(\mathrm{LL})}$ | $\mathrm{CM}_n^{(\mathrm{LLS})}$ | $\mathrm{KS}_n^{(\mathrm{E})}$ | $\mathrm{KS}_n^{(\mathrm{LL})}$ | $\mathrm{KS}_n^{(\mathrm{LLS})}$ | $Q_n^{(\mathrm{E})}$ | $Q_n^{(\mathrm{LL})}$ | $Q_n^{(\mathrm{LLS})}$ |
| Clayton | Clayton | *4.9* | *4.3* | *4.9* | *5.1* | *3.8* | *4.5* | *5.0* | *4.9* | *5.6* |
| | Gumbel | 86.7 | 91.3 | **93.3** | 61.1 | 81.8 | **84.6** | 84.1 | 91.7 | **93.6** |
| | Frank | 54.4 | **63.5** | 61.4 | 33.5 | **64.1** | 61.2 | 51.6 | **62.2** | 59.7 |
| | Plackett | 56.3 | 61.0 | **63.5** | 34.2 | **57.5** | 57.2 | 53.0 | 61.8 | **62.7** |
| | Normal | 49.9 | 59.0 | **61.6** | 28.2 | **54.8** | 52.8 | 46.6 | 59.6 | **62.5** |
| | Student, 4 df | 55.9 | 65.5 | **67.0** | 34.6 | **40.8** | 40.7 | 52.1 | 68.0 | **69.5** |
| Gumbel | Clayton | 72.2 | 85.2 | **86.7** | 48.3 | 77.5 | **78.1** | 74.0 | 83.3 | **85.6** |
| | Gumbel | *4.9* | *5.1* | *5.6* | *5.0* | *4.9* | *4.7* | *4.5* | *4.9* | *5.4* |
| | Frank | 14.7 | 18.2 | **20.1** | 11.0 | 18.2 | **22.8** | **18.4** | 15.7 | 17.3 |
| | Plackett | 13.7 | 16.7 | **17.4** | 9.5 | 17.7 | **17.9** | 16.2 | 15.7 | **16.3** |
| | Normal | 10.0 | **15.6** | 15.5 | 8.2 | **16.8** | 16.1 | 12.6 | **15.5** | 14.1 |
| | Student, 4 df | 12.8 | 21.5 | **22.6** | 7.3 | 11.3 | **13.3** | 14.3 | 23.3 | **24.5** |
| Frank | Clayton | 41.3 | 49.9 | **52.1** | 25.5 | 37.6 | **39.1** | 42.1 | 48.2 | **50.8** |
| | Gumbel | 31.8 | 47.5 | **50.1** | 18.1 | **28.5** | 28.3 | 25.4 | 45.7 | **47.5** |
| | Frank | *4.7* | *4.9* | *4.9* | *5.1* | *4.9* | *4.5* | *4.6* | *4.6* | *5.1* |
| | Plackett | 5.7 | **6.9** | 6.2 | 4.7 | 5.0 | **5.8** | 5.4 | 6.7 | **7.0** |
| | Normal | 8.5 | **13.4** | 12.6 | 9.7 | **12.1** | 11.5 | 7.6 | **14.4** | 12.8 |
| | Student, 4 df | 18.1 | **34.5** | 33.9 | 10.1 | **18.1** | 17.1 | 15.4 | **36.6** | 35.2 |
| Normal | Clayton | 34.3 | 36.9 | **40.1** | 18.3 | 25.6 | **28.7** | 36.2 | 36.0 | **38.9** |
| | Gumbel | 26.1 | **35.1** | 33.7 | 12.0 | **20.0** | 18.7 | 21.0 | **34.0** | 32.6 |
| | Frank | **7.6** | 3.7 | 3.9 | **7.1** | 3.4 | 4.3 | **7.9** | 4.2 | 4.3 |
| | Plackett | **8.1** | 4.5 | 5.4 | **8.0** | 3.6 | 5.1 | **8.6** | 5.7 | 6.1 |
| | Normal | *4.8* | *5.2* | *5.5* | *5.0* | *5.2* | *4.5* | *5.8* | *4.9* | *5.1* |
| | Student, 4 df | 11.7 | **14.8** | 14.5 | 6.5 | **9.2** | 9.0 | 9.5 | 17.9 | **19.1** |
| Student | Clayton | 29.0 | 29.5 | **34.1** | 22.1 | 33.2 | **36.3** | **32.8** | 26.2 | 30.2 |
| | Gumbel | 20.7 | 26.1 | **26.6** | 12.1 | 22.7 | **25.5** | 17.9 | 17.9 | **22.7** |
| | Frank | **9.0** | 7.1 | 6.1 | 9.1 | 8.0 | **9.2** | **9.7** | 3.6 | 3.4 |
| | Plackett | **7.6** | 6.4 | 4.5 | **8.6** | 6.6 | 7.7 | **7.6** | 3.9 | 2.9 |
| | Normal | **4.7** | 4.3 | 3.5 | 6.4 | **6.9** | 6.7 | **4.9** | 3.9 | 2.9 |
| | Student, 4 df | *4.7* | *4.7* | *4.9* | *5.3* | *4.8* | *4.7* | *4.4* | *4.7* | *5.4* |

where $F_n^*(u) = \frac{1}{n} \sum_{i=1}^n \mathbb{I}\{U_i \le u\}$ and $G_n^*(v) = \frac{1}{n} \sum_{i=1}^n \mathbb{I}\{V_i \le v\}$. The convergence of the smoothed copula process $\mathbb{C}_n^{(\mathrm{LL})} = \sqrt{n}(\hat{C}_n^{(\mathrm{LL})} - C)$ to a Gaussian process given by (17) will now follow by the results of Fermanian, Radulovič and Wegkamp (2004).

STEP 1. We consider the following class of functions from $[0,1]^2$ to $[0,1]$:
$$\mathcal{F} = \left\{ (w_1, w_2) \mapsto K_{u,h}\left(\frac{u - \zeta_1(w_1)}{h}\right) K_{v,h}\left(\frac{v - \zeta_2(w_2)}{h}\right), \right.$$

TABLE 3
*Percentage of rejection of $H_0$ by various tests for samples of size $n = 150$ arising from different copula models with $\tau = 0.50$*

| Copula under $H_0$ | True copula | Cramér–von Mises | | | Kolmogorov–Smirnov | | | MISE | | |
|---|---|---|---|---|---|---|---|---|---|---|
| | | $\text{CM}_n^{(E)}$ | $\text{CM}_n^{(LL)}$ | $\text{CM}_n^{(LLS)}$ | $\text{KS}_n^{(E)}$ | $\text{KS}_n^{(LL)}$ | $\text{KS}_n^{(LLS)}$ | $Q_n^{(E)}$ | $Q_n^{(LL)}$ | $Q_n^{(LLS)}$ |
| Clayton | Clayton | *5.3* | *5.2* | *5.4* | *5.5* | *5.3* | *5.7* | *4.8* | *5.6* | *5.8* |
| | Gumbel | 99.9 | 100.0 | 99.9 | 98.9 | 99.5 | **99.9** | 100.0 | 100.0 | 99.9 |
| | Frank | 95.9 | **96.4** | 96.2 | 82.5 | **98.3** | 95.6 | 91.1 | **93.3** | 92.3 |
| | Plackett | 95.6 | **96.7** | 96.0 | 75.3 | **94.7** | 89.7 | 92.3 | **95.7** | 94.5 |
| | Normal | 94.4 | 96.3 | **96.9** | 75.0 | **91.9** | 91.1 | 89.9 | 94.5 | **94.9** |
| | Student, 4 df | 94.9 | 96.7 | **97.2** | 77.9 | 86.4 | **88.2** | 92.7 | 96.2 | **96.3** |
| Gumbel | Clayton | 99.6 | 99.7 | 99.7 | 94.3 | 98.7 | **98.8** | **98.9** | 99.3 | 99.1 |
| | Gumbel | *4.5* | *4.7* | **4.9** | *5.0* | *4.9* | *5.4* | *4.8* | *3.9* | *4.6* |
| | Frank | 40.5 | **48.2** | 39.7 | 29.6 | **48.1** | 40.4 | **41.1** | 35.3 | 30.6 |
| | Plackett | 29.4 | **33.3** | 30.7 | 18.9 | **26.5** | 22.8 | **31.0** | 30.0 | 27.6 |
| | Normal | 18.8 | **26.4** | 25.1 | 14.6 | **25.3** | 23.8 | **22.3** | 22.1 | 21.8 |
| | Student, 4 df | 22.3 | 27.8 | **29.2** | 11.7 | **19.1** | 16.5 | 23.6 | 28.5 | **29.0** |
| Frank | Clayton | 89.6 | 88.9 | **91.9** | 68.1 | 72.7 | **75.5** | 85.1 | **86.3** | 85.5 |
| | Gumbel | 63.8 | 71.0 | **74.1** | 39.3 | 44.6 | **47.3** | 50.5 | **68.9** | 65.7 |
| | Frank | *5.3* | *4.9* | *5.2* | *5.1* | *5.2* | *5.0* | *5.1* | *4.9* | *5.0* |
| | Plackett | 8.4 | 10.4 | **12.1** | 5.4 | **6.9** | 6.7 | 8.3 | **15.2** | 8.5 |
| | Normal | 19.6 | 26.0 | **29.5** | 17.6 | **26.9** | 25.5 | 16.5 | **25.7** | 17.5 |
| | Student, 4 df | 35.0 | 44.9 | **52.8** | 17.9 | 27.8 | **29.4** | 29.0 | **51.0** | 46.2 |
| Normal | Clayton | **83.0** | 78.3 | 82.9 | 55.8 | **66.6** | 66.1 | **79.6** | 76.1 | 79.4 |
| | Gumbel | 41.7 | 39.5 | **44.3** | 18.3 | 22.7 | **26.6** | 32.4 | 39.6 | **41.2** |
| | Frank | **21.2** | 20.1 | 14.9 | **15.1** | 11.8 | 11.0 | **19.3** | 14.4 | 10.5 |
| | Plackett | **12.0** | 7.4 | 7.8 | **7.7** | 5.8 | 4.5 | **13.4** | 11.9 | 9.9 |
| | Normal | *4.8* | *5.1* | *5.5* | *4.7* | *5.0* | *5.4* | *5.8* | *4.5* | *4.6* |
| | Student, 4 df | 8.1 | 6.4 | **8.3** | 4.0 | 4.5 | **4.8** | 7.6 | 11.5 | **12.5** |
| Student | Clayton | 80.6 | 78.4 | **83.7** | 62.1 | 74.9 | **76.5** | **79.0** | 72.3 | 75.0 |
| | Gumbel | 36.5 | 39.1 | **39.6** | 20.5 | 31.3 | **32.4** | 25.7 | 25.7 | **31.0** |
| | Frank | 28.5 | **30.5** | 23.1 | 18.8 | 18.0 | **20.8** | **25.0** | 15.0 | 11.3 |
| | Plackett | **13.5** | **13.5** | 8.9 | **10.0** | 6.9 | 7.3 | **11.0** | 6.8 | 6.6 |
| | Normal | 5.1 | **6.4** | 5.5 | 7.6 | **7.8** | 7.5 | **5.1** | 3.4 | 4.0 |
| | Student, 4 df | *4.7* | *4.9* | *5.1* | *4.9* | *5.0* | *4.9* | *5.2* | *4.1* | *4.7* |

(27)
$$(u,v) \in [0,1]^2, h \in \left[0, \frac{1}{4}\right], \zeta_1, \zeta_2 : [0,1] \to [0,1] \text{ nondecreasing} \bigg\}.$$

As each function $f$ from $\mathcal{F}$ is characterized by a quintuple $(u, v, h, \zeta_1, \zeta_2)$, the empirical process indexed by $\mathcal{F}$ can be written as

$$Z_n(f) = Z_n(u, v, h, \zeta_1, \zeta_2)$$

TABLE 4
*Percentage of rejection of $H_0$ by various tests for samples of size $n = 150$ arising from different copula models with $\tau = 0.75$*

| Copula under $H_0$ | True copula | Cramér–von Mises | | | Kolmogorov–Smirnov | | | MISE | | |
|---|---|---|---|---|---|---|---|---|---|---|
| | | $\mathrm{CM}_n^{(\mathrm{E})}$ | $\mathrm{CM}_n^{(\mathrm{LL})}$ | $\mathrm{CM}_n^{(\mathrm{LLS})}$ | $\mathrm{KS}_n^{(\mathrm{E})}$ | $\mathrm{KS}_n^{(\mathrm{LL})}$ | $\mathrm{KS}_n^{(\mathrm{LLS})}$ | $Q_n^{(\mathrm{E})}$ | $Q_n^{(\mathrm{LL})}$ | $Q_n^{(\mathrm{LLS})}$ |
| Clayton | Clayton | *5.3* | *5.6* | *5.3* | *5.1* | *4.9* | *4.6* | *3.3* | *4.6* | *4.0* |
| | Gumbel | 100.0 | 100.0 | 100.0 | 99.9 | 100.0 | 100.0 | 100.0 | 100.0 | 100.0 |
| | Frank | **98.8** | 98.3 | 98.5 | 83.7 | **98.1** | 97.2 | 94.0 | **95.4** | 95.1 |
| | Plackett | 99.5 | 99.2 | **99.9** | 85.1 | 90.6 | **93.9** | 97.0 | 98.3 | **99.1** |
| | Normal | 99.8 | 99.5 | **99.9** | 90.5 | 93.6 | **96.8** | 98.0 | 98.4 | **99.2** |
| | Student, 4 df | 99.9 | 99.7 | **100.0** | 92.7 | 91.4 | **97.7** | 98.5 | **99.5** | 99.5 |
| Gumbel | Clayton | **99.9** | 99.5 | **99.9** | 95.8 | **98.5** | **98.5** | **99.6** | 99.1 | 99.0 |
| | Gumbel | *4.5* | *4.6* | *4.8* | *4.7* | *4.6* | *4.9* | *5.0* | *3.2* | *3.7* |
| | Frank | 53.3 | **54.5** | 47.0 | 25.6 | 38.1 | **38.3** | **40.7** | 26.6 | 23.7 |
| | Plackett | **24.3** | 23.5 | 19.1 | 6.6 | **9.3** | **9.3** | **29.5** | 24.9 | 25.1 |
| | Normal | 12.4 | 13.1 | **13.9** | 11.3 | **13.6** | 13.3 | **12.3** | 8.3 | 7.7 |
| | Student, 4 df | 15.6 | 15.7 | **20.1** | 8.5 | 10.4 | **10.8** | **16.2** | 15.1 | 14.9 |
| Frank | Clayton | 96.7 | 91.3 | **96.9** | 57.6 | 63.7 | **64.9** | **90.8** | 86.2 | 89.0 |
| | Gumbel | 81.6 | 80.5 | **87.8** | 36.9 | **41.7** | 39.7 | 61.1 | 73.7 | **75.3** |
| | Frank | *4.4* | *4.6* | *4.4* | *4.7* | *4.5* | *4.8* | *4.7* | *3.0* | *3.7* |
| | Plackett | 19.6 | 20.7 | **27.5** | 5.9 | **8.7** | 7.5 | 34.5 | 45.7 | **48.3** |
| | Normal | 40.7 | 41.1 | **52.7** | 28.5 | **33.3** | 30.2 | 31.9 | 38.7 | **42.9** |
| | Student, 4 df | 58.4 | 57.9 | **72.6** | 26.7 | **33.6** | 30.3 | 50.6 | 63.1 | **64.6** |
| Normal | Clayton | **93.4** | 88.9 | 90.5 | 66.7 | 77.9 | **79.7** | **88.4** | 83.3 | 82.6 |
| | Gumbel | 41.1 | 38.3 | **43.4** | 13.0 | 19.6 | **20.9** | 24.3 | 33.5 | **35.0** |
| | Frank | 46.1 | **46.7** | 37.7 | 17.5 | 18.2 | **22.9** | **31.0** | 24.1 | 18.4 |
| | Plackett | **15.2** | 11.6 | 9.2 | 3.1 | 3.1 | **4.0** | 24.4 | **27.7** | 23.7 |
| | Normal | *4.7* | *4.8* | *4.4* | *4.4* | *4.6* | *4.7* | *5.2* | *3.5* | *3.5* |
| | Student, 4 df | **6.9** | 6.3 | 6.8 | **4.7** | 3.9 | 4.2 | 7.0 | **10.1** | 9.6 |
| Student | Clayton | **92.8** | 89.4 | 89.9 | 73.7 | 84.4 | **86.7** | **85.8** | 75.6 | 74.8 |
| | Gumbel | **37.3** | 34.5 | 37.0 | 17.3 | 26.3 | **26.9** | 18.4 | 18.4 | **21.6** |
| | Frank | **52.2** | 51.8 | 45.5 | 24.5 | 28.1 | **33.6** | **30.4** | 20.7 | 14.5 |
| | Plackett | **16.4** | 15.7 | 10.3 | 4.3 | 3.9 | **5.8** | 5.8 | **12.7** | 11.6 |
| | Normal | 4.5 | **4.7** | 3.5 | 6.3 | 6.9 | **7.8** | **2.5** | 2.1 | 1.9 |
| | Student, 4 df | *4.3* | *4.9* | *4.4* | *4.7* | *5.1* | *4.9* | *4.8* | *3.6* | *3.1* |

$$= \frac{1}{\sqrt{n}} \sum_{i=1}^{n} K_{u,h}\left(\frac{u - \zeta_1(U_i)}{h}\right) K_{v,h}\left(\frac{v - \zeta_2(V_i)}{h}\right).$$

Put $\bar{Z}_n = Z_n - \mathrm{E} Z_n$ and note that

(28)
$$A_n^{h_n}(u,v) = \bar{Z}_n(f_1^n) - \bar{Z}_n(f_2),$$
where $f_1^n = (u, v, h_n, F_n(F^{-1}), G_n(F^{-1})), \ f_2 = (u, v, 0, \mathbb{I}, \mathbb{I})$

with $\mathbb{I}$ being the identity function on the interval $[0,1]$.

Lemma A.1, which is given below, states that the set of functions $\mathcal{F}$ is Donsker. Indeed, $\mathcal{F}$ is a subset of $\mathcal{F}^*$ in Lemma A.1, taking $b(\cdot) = 1$ and $u_0 = u$ and $v_0 = v$. This implies the weak convergence of the process $\bar{Z}_n(f), f \in \mathcal{F}$, which further implies that the process $\bar{Z}_n$ is asymptotically uniformly $\rho$-equicontinuous in probability [see pages 37–41 of van der Vaart and Wellner (1996)] with semimetric $\rho$ given by

$$\rho(f, f') = \mathrm{E}\bigg[K_{u,h}\bigg(\frac{u - \zeta_1(U)}{h}\bigg)K_{v,h}\bigg(\frac{v - \zeta_2(V)}{h}\bigg)$$
$$- K_{u',h'}\bigg(\frac{u' - \zeta_1'(U)}{h'}\bigg)K_{v',h'}\bigg(\frac{v' - \zeta_2'(V)}{h'}\bigg)\bigg]^2.$$

Using this asymptotic uniform $\rho$-equicontinuity and (28), we get that $\sup_{u,v} |A_n^{h_n}| = o_p(1)$, provided that $\sup_{u,v} \rho(f_1^n, f_2)$ converges to zero in probability, where $f_1^n$ and $f_2$ are given in (28) [for details consult the proof in van der Vaart (1994)].

Put $M = \sup_{u,h,x} |k_{u,h}(x)|$, where $k_{u,h}$ is defined in (7), and denote

$$A_\varepsilon = \{|\hat{U} - U| > \varepsilon \text{ or } |\hat{V} - V| > \varepsilon\}.$$

The consistency of $\hat{F}_n$ and $\hat{G}_n$ yields that, for every $\varepsilon > 0$, for all sufficiently large $n$,

$$P\Big[\max\Big\{\sup_{x \in \mathbb{R}} |\hat{F}_n(x) - F(x)|, \sup_{y \in \mathbb{R}} |\hat{G}_n(y) - G(y)|\Big\} > \varepsilon\Big] < \varepsilon,$$

which further implies that, for all sufficiently large $n$,

(29) $\quad P(A_\varepsilon) = P[\max\{|\hat{F}_n(X) - F(X)|, |\hat{G}_n(Y) - G(Y)|\} > \varepsilon] < \varepsilon.$

Now, we can bound

$$\rho(f_1^n, f_2) = \mathrm{E}\bigg[K_{u,h_n}\bigg(\frac{u - \hat{U}}{h_n}\bigg)K_{v,h_n}\bigg(\frac{v - \hat{V}}{h_n}\bigg) - \mathbb{I}\{U \leq u, V \leq v\}\bigg]^2$$
$$\leq \mathrm{E}\bigg[K_{u,h}\bigg(\frac{u - \hat{U}}{h}\bigg)K_{v,h}\bigg(\frac{v - \hat{V}}{h}\bigg) - \mathbb{I}\{U \leq u, V \leq v\}\bigg]^2 \mathbb{I}_{A_\varepsilon^c} + \mathbb{I}_{A_\varepsilon}$$
$$\leq M^4 \mathrm{E}[|\mathbb{I}\{\hat{U} \leq u - h_n\} - \mathbb{I}\{U \leq u\}|$$
$$\qquad + |\mathbb{I}\{\hat{U} \leq u + h_n\} - \mathbb{I}\{U \leq u\}|$$
$$\qquad + |\mathbb{I}\{\hat{V} \leq v - h_n\} - \mathbb{I}\{V \leq v\}|$$
$$\qquad + |\mathbb{I}\{\hat{V} \leq v + h_n\} - \mathbb{I}\{V \leq v\}|]\mathbb{I}_{A_\varepsilon^c} + \mathbb{I}_{A_\varepsilon}$$
$$\leq 4M^4(\varepsilon + h_n) + \mathbb{I}_{A_\varepsilon}.$$

As the above bound holds uniformly in $(u, v)$ and by (29), for all sufficiently large $n$, we have $P(A_\varepsilon) < \varepsilon$. Since $\varepsilon$ can be arbitrarily small, this implies $\sup_{u,v} \rho(f_1^n, f_2) = o_p(1)$.

LEMMA A.1. *Suppose that the function $k$ is of bounded variation and $\int k(x)\,dx = 1$. Then, the set of functions from $[0,1]^2$ to $[0,1]$*

$$\mathcal{F}^* = \bigg\{(w_1, w_2) \mapsto K_{u_0, h}\bigg(\frac{u - \zeta_1(w_1)}{b(u_0)h}\bigg) K_{v_0, h}\bigg(\frac{v - \zeta_2(w_2)}{b(v_0)h}\bigg),$$

$$(u_0, v_0), (u, v) \in [0, 1]^2,$$

$$h \in \bigg[0, \frac{1}{4}\bigg], \zeta_1, \zeta_2 : [0, 1] \to [0, 1] \text{ nondecreasing}\bigg\},$$

*where $b(w) = 1$ or $b(w) = \min\{\sqrt{w}, \sqrt{1-w}\}$, is Donsker.*

*Consequently, the family $\mathcal{F}$ in (27) is Donsker.*

PROOF. Note that the class of functions

$$\mathcal{G}_1 = \{(w_1, w_2) \mapsto \mathbb{I}\{\zeta_1(w_1) \le a, \zeta_2(w_2) \le b\},$$

$$a, b \in \mathbb{R}, \zeta_1, \zeta_2 : [0, 1] \to [0, 1] \text{ nondecreasing}\}$$

is a subset of the class of indicators

$$\mathcal{G}_2 = \{(w_1, w_2) \mapsto \mathbb{I}\{w_1 < (\le)a, w_2 < (\le)b\}, a, b \in \mathbb{R}\}.$$

But this implies that $\mathcal{G}_1$ is a Donsker class [see Example 2.5.4 of van der Vaart and Wellner (1996)].

As the set $\mathcal{G}_1$ is closed under translation, we know by the beginning of the proof of van der Vaart (1994) that the set of functions

$$(30) \qquad \mathcal{H} = \bigg\{\int f(x+y)\,d\mu(y), f \in \mathcal{G}_1, \mu \in \mathcal{M}_B\bigg\}$$

is a Donsker class, where $\mathcal{M}_B$ is a family of all signed measures (on $\mathbb{R}^2$) of total mass bounded by a fixed constant $B$.

Let us introduce the set of signed measures

$$\mathcal{M}_0 = \bigg\{(-\infty, w_1] \times (-\infty, w_2] \mapsto K_{u_0, h}\bigg(\frac{w_1}{b(u_0)h}\bigg) K_{v_0, h}\bigg(\frac{w_2}{b(v_0)h}\bigg),$$

$$(u_0, v_0) \in [0, 1]^2, h \in \bigg[0, \frac{1}{4}\bigg]\bigg\}.$$

If $k$ is of bounded variation, then, by taking sufficiently large $B$, we ensure that $\mathcal{M}_0 \subset \mathcal{M}_B$. Further, if $x$ stands for $(w_1, w_2)$ and $y$ for $(y_1, y_2)$, then, for $f \in \mathcal{G}_1$ and $\mu \in \mathcal{M}_0$, we get

$$\int f(x+y)\,d\mu(y)$$

$$= \int\int \mathbb{I}\{\zeta_1(w_1) + y_1 \le u,$$

IMPROVED KERNEL ESTIMATION OF COPULAS 25

$$\zeta_2(w_2) + y_2 \leq v\} d\bigg(K_{u_0,h}\bigg(\frac{y_1}{b(u_0)h}\bigg)K_{v_0,h}\bigg(\frac{y_2}{b(v_0)h}\bigg)\bigg)$$

$$= K_{u_0,h}\bigg(\frac{u - \zeta_1(w_1)}{b(u_0)h}\bigg)K_{v_0,h}\bigg(\frac{v - \zeta_2(w_2)}{b(v_0)h}\bigg),$$

which is a Donsker class. As $\mathcal{F}^*$ includes $\mathcal{F}$ of (27) (consider $u_0 = u$ and $v_0 = v$), the family $\mathcal{F}$ is a Donsker class as well. $\square$

STEP 2. Now, we can turn our attention to the process $C_n^{h_n}$ given by (26). For $(u,v) \in \mathbb{R}^2$, define

$$C^*_{\tilde{F},\tilde{G}}(u,v) = C(F(\tilde{F}^{-1}(u^*)), G(\tilde{G}^{-1}(v^*))),$$

where $w^* = \max\{\min\{w,1\},0\}$. Note that

$$\mathbb{E}\mathbb{I}\{\hat{U} \leq u, \hat{V} \leq v\} = \mathbb{E}\mathbb{I}\{U \leq F(\hat{F}_n^{-1}(u)), \hat{V} \leq G(\hat{G}_n^{-1}(u))\} + O(n^{-1})$$
$$= C^*_{\hat{F}_n,\hat{G}_n}(u,v) + O(n^{-1}),$$

where the remainder term $O(n^{-1})$ disappears if we do some smoothing on the first stage; that is, if $b_{1n}, b_{2n} > 0$. As

$$\mathbb{E}K_{u,h}\bigg(\frac{u - \hat{U}}{h}\bigg)K_{v,h}\bigg(\frac{v - \hat{V}}{h}\bigg)$$

$$= \mathbb{E}\int_{-1}^{1}\int_{-1}^{1}\mathbb{I}\{\hat{U}_i \leq u - th_n, \hat{V}_i \leq u - sh_n\}k_{u,h}(s)k_{v,h}(t)\,dt\,ds$$

$$= \int_{-1}^{1}\int_{-1}^{1}C^*_{\hat{F}_n,\hat{G}_n}(u - th_n, v - sh_n)k_{u,h}(s)k_{v,h}(t)\,dt\,ds + O(n^{-1}),$$

it will be useful to have a closer look at the process $\{C^*_{\hat{F}_n,\hat{G}_n}(u,v) \in [0,1]^2\}$. In the following, we will prove that, uniformly in $(u,v)$,

$$\sqrt{n}(C^*_{\hat{F}_n,\hat{G}_n}(u,v) - C(u,v))$$

(31)
$$= -C_u(u,v)\frac{1}{\sqrt{n}}\sum_{i=1}^{n}[\mathbb{I}\{U_i \leq u\} - u]$$

$$\quad - C_v(u,v)\frac{1}{\sqrt{n}}\sum_{i=1}^{n}[\mathbb{I}\{V_i \leq v\} - v] + o_p(1).$$

Let $D[a,b]$ be the Banach space of all cadlag functions on an interval $[a,b]$ equipped with the uniform norm.

LEMMA A.2. *Let $F$ be a continuous distribution function. Then, the map $\tilde{F} \mapsto F \circ \tilde{F}^{-1}$ as a map $D[0,1] \mapsto \ell^\infty[0,1]$ is Hadamard-differentiable at $\tilde{F} = F$ tangentially to the set of functions*

$$\alpha \in \mathrm{E}_F = \{\alpha(x) = \beta(F(x)), x \in \mathbb{R}, \beta \in C[0,1], \beta(0) = \beta(1) = 0\}.$$

*The derivative is given by* $-\alpha \circ F^{-1}$.

PROOF. Let $\alpha_t$ converge uniformly to $\alpha \in \mathrm{E}_F$. Put $F_t = F + t\alpha_t$. As

$$\sup_{u \in [0,1]} |F_t(F_t^{-1}(u)) - u| = \sup_{u \in [0,1]} |F_t(u+) - F_t(u-)|,$$

the function $F$ is continuous and the function $\alpha_t$ converges uniformly to a bounded and continuous function, we get that $F_t(F_t^{-1}(u)) - u = o(t)$ uniformly in $u$.

Thus, we can calculate

$$\frac{1}{t}[F(F_t^{-1})(u) - F(F^{-1}(u)) + t\alpha(F^{-1}(u))]$$
$$= \alpha(F^{-1}(u)) - \alpha_t(F_t^{-1}(u)) + o(1)$$
$$= [\alpha(F^{-1}(u)) - \alpha(F_t^{-1}(u))] + [\alpha(F_t^{-1}(u)) - \alpha_t(F_t^{-1}(u))] + o(1).$$

As $\alpha_t \to \alpha$ uniformly, the second term converges to zero uniformly in $u$. By using the representation $\alpha(x) = \beta(F(x))$, we see that, to ensure a uniform convergence of the first term to zero, we need to show that $F(F_t^{-1}(u)) \to u$ uniformly. But this follows by a simple calculation, which yields

$$|F(F_t^{-1}(u)) - F(F^{-1}(u))|$$
$$= |t\alpha_t(F_t^{-1}(u))| + o(t)$$
$$\leq |t||\alpha_t(F_t^{-1}(u)) - \alpha(F_t^{-1}(u))| + |t||\alpha(F_t^{-1}(u))| + o(t) = O(t),$$

uniformly in $u$. $\square$

REMARK. As (16) implies that $\sqrt{n}(\hat{F}_n - F)$ converges in distribution to a Gaussian process $B \circ F$, where $B$ is a standard Brownian motion on the interval $[0, 1]$, the Hadamard-differentiability tangentially to $\mathrm{E}_F$ given in Lemma A.2 is exactly what is needed to derive asymptotic distribution of the process $\sqrt{n}(F(F_n^{-1}(u)) - u)$.

Similarly, we can prove that the mapping $\tilde{G} \mapsto G \circ \tilde{G}^{-1}$ is Hadamard-differentiable at $\tilde{G} = G$ tangentially to the set of functions

$$\mathrm{E}_G = \{\alpha(x) = \beta(G(x)), x \in \mathbb{R}, \beta \in C[0,1], \beta(0) = \beta(1) = 0\}.$$

The proof of the following lemma follows easily by applying Lemma A.2 and the chain rule.

LEMMA A.3. *Let the copula $C$ have continuous partial derivatives on $[a,b] \times [c,d] \subset [0,1]^2$, and $F, G$ are continuous; then, the map $(\tilde{F}, \tilde{G}) \mapsto C^*_{\tilde{F},\tilde{G}}$ as a map $D[a,b] \times D[a,b] \mapsto \ell^\infty([a,b] \times [c,d])$ is Hadamard-differentiable at* the point $(\tilde{F}, \tilde{G}) = (F, G)$ tangentially to the set of functions $(\alpha_1, \alpha_2) \in \mathrm{E}_F \times \mathrm{E}_G$. The derivative is given by

(32) $\quad \phi'(\alpha_1, \alpha_2) = -C_u \circ \alpha_1 \circ F^{-1} - C_v \circ \alpha_2 \circ G^{-1} = -C_u \circ \beta_1 - C_v \circ \beta_2,$

where $\beta_1 = \alpha_1 \circ F^{-1}$ and $\beta_2 = \alpha_2 \circ G^{-1}$.

Lemma A.3, together with Theorem 3.9.4 of van der Vaart and Wellner (1996), imply that representation (31) holds uniformly for $(u, v) \in [a, b] \times [c, d]$. Unfortunately, many of the most popular families (e.g., Clayton, Gumbel, normal) do not have continuous $C_u$ and $C_v$ at some of the points $\{(0,0), (0,1), (1,0), (1,1)\}$. The following lemma takes care of this situation.

LEMMA A.4. *Let the distribution $H$ have continuous margins $F, G$ and a copula function whose first derivatives are continuous on $[0,1]^2 \setminus \{(0,0), (0,1), (1,0), (1,1)\}$. Then, representation (31) holds uniformly in $(u,v) \in [0,1]^2$.*

PROOF. Suppose, for simplicity, that the point of discontinuity is only at $(0,0)$ [other points $(0,1), (1,0), (1,1)$ might be handled in a similar way]. Let us denote

$$Z_n(u,v) = \sqrt{n}(C^*_{\hat{F}_n, \hat{G}_n}(u,v) - C(u,v))$$
(33)
$$+ C_u(u,v) \frac{1}{\sqrt{n}} \sum_{i=1}^n [\mathbb{I}\{U_i \leq u\} - u]$$
$$+ C_v(u,v) \frac{1}{\sqrt{n}} \sum_{i=1}^n [\mathbb{I}\{V_i \leq v\} - v].$$

Let $\varepsilon > 0$ be given. As all the process

$$X_n^1(u) = \sqrt{n}[F(\hat{F}_n^{-1}(u)) - u], \qquad X_n^3(u) = \frac{1}{\sqrt{n}} \sum_{i=1}^n [\mathbb{I}\{U_i \leq u\} - u],$$

$$X_n^2(u) = \sqrt{n}[G(\hat{G}_n^{-1}(u)) - u], \qquad X_n^4(u) = \frac{1}{\sqrt{n}} \sum_{i=1}^n [\mathbb{I}\{V_i \leq u\} - u]$$

converge to a Brownian motion, we can find $\delta_\varepsilon$ and $n_\varepsilon$ such that, for all $n > n_\varepsilon$,

$$P\left(\sup_{u \leq \delta_\varepsilon} |X_n^j(u)| \geq \frac{\varepsilon}{4}\right) < \frac{\varepsilon}{4}, \qquad j = 1, \ldots, 4.$$

As $C_u, C_v$ are bounded by 1, the triangular inequality implies that, for all $n > n_\varepsilon$,

$$P\left(\sup_{u,v \leq \delta_\varepsilon} |Z_n(u,v)| \geq \varepsilon\right) \leq \sum_{j=1}^4 P\left(\sup_{u \leq \delta_\varepsilon} |X_n^j(u)| \geq \frac{\varepsilon}{4}\right) < \varepsilon.$$

Next, the existence of $n'_\varepsilon$ such that, for all $n > n'_\varepsilon$,

$$P\Big(\sup_{u,v \in A_\varepsilon} |Z_n(u,v)| \geq \varepsilon\Big) < \varepsilon \quad \text{with } A_\varepsilon = [0,1]^2 \setminus [0,\delta_\varepsilon]^2,$$

follows by Lemma A.3 applied to rectangles $[0,\delta_\varepsilon] \times [\delta_\varepsilon, 1]$ and $[\delta_\varepsilon, 1] \times [0,1]$. Thus, for $n > \max\{n_\varepsilon, n'_\varepsilon\} : P(\sup_{u,v} |Z_n(u,v)| \geq \varepsilon) < \varepsilon$, which proves the lemma. $\square$

Combining (31), Lemma A.4, the fact that $h_n \to 0$ and asymptotic equicontinuity of the processes $\mathbb{U}_n(u) = \frac{1}{\sqrt{n}} \sum_{i=1}^n [\mathbb{I}\{U_i \leq u\} - u]$, $\mathbb{V}_n(v) = \frac{1}{\sqrt{n}} \sum_{i=1}^n [\mathbb{I}\{V_i \leq v\} - v]$ yields

$$\begin{aligned}
C_n^{h_n}(u,v) &= \sqrt{n} \mathbb{E}\Big[K_{u,h_n}\Big(\frac{u-\hat{U}}{h_n}\Big) K_{v,h_n}\Big(\frac{v-\hat{V}}{h_n}\Big) - C(u,v)\Big] \\
&= -C_u(u,v) \mathbb{U}_n(u) - C_v(u,v) \mathbb{V}_n(v) + \sqrt{n} D_n(u,v) + o_P(1),
\end{aligned} \tag{34}$$

where the bias term $D_n$ is given by

$$\begin{aligned}
D_n(u,v) = \int_{-1}^1 \int_{-1}^1 \sqrt{n}[C(u-th_n, v-sh_n) - C(u,v)] \\
\times k_{u,h}(s) k_{v,h}(t)\, dt\, ds.
\end{aligned} \tag{35}$$

If copula $C$ has bounded second order partial derivatives on $[0,1]^2$, then

$$\sqrt{n} \sup_{u,v} |D_n(u,v)| = O(n^{1/2} h_n^2) = o(1). \tag{36}$$

Finally, combining (23), (25), (26), (34) and (36) yields

$$\begin{aligned}
\sqrt{n}(\hat{C}_n^{(\text{LL})}(u,v) - C(u,v)) \\
= \frac{1}{\sqrt{n}} \sum_{i=1}^n [\mathbb{I}\{U_i \leq u, V_i \leq v\} - C(u,v)] \\
- C_u(u,v) \frac{1}{\sqrt{n}} \sum_{i=1}^n [\mathbb{I}\{U_i \leq u\} - u] \\
- C_v(u,v) \frac{1}{\sqrt{n}} \sum_{i=1}^n [\mathbb{I}\{V_i \leq v\} - v] + o_P(1).
\end{aligned} \tag{37}$$

**A.2. Weak convergence of the process $\mathbb{C}_n^{(\text{MR})}$.** Here, we adapt the foregoing proof for the mirrored-type kernel estimator $\hat{C}_n^{(\text{MR})}$ given in (13).

STEP 1. At first, we rewrite

$$(38) \quad \hat{C}_n^{(\mathrm{MR})} = \sum_{\ell=1}^{9}[Z_n(\ell,u,v) - Z_n(\ell,u,0) - Z_n(\ell,0,v) + Z_n(\ell,0,0)],$$

where

$$Z_n(\ell,u,v) = \frac{1}{n}\sum_{i=1}^{n} K\left(\frac{u - \hat{U}_i^{(\ell)}}{h_n}\right) K\left(\frac{v - \hat{V}_i^{(\ell)}}{h_n}\right).$$

Let us define

$$Z_{0n}(\ell,u,v) = \frac{1}{n}\sum_{i=1}^{n} \mathbb{I}\{U_i^{(\ell)} \leq u, V_i^{(\ell)} \leq v\}.$$

Similarly as in Step 1 of the proof of Appendix A.1 [weak convergence of $\mathbb{C}_n^{(\mathrm{LL})}$], we can show that, for each $\ell = 1,\ldots,9$,

$$(39) \quad \sup_{u,v}|\sqrt{n}(Z_n(\ell,u,v) - \mathrm{E}Z_n(\ell,u,v)) - \sqrt{n}(Z_{0n}(\ell,u,v) - \mathrm{E}Z_{0n}(\ell,u,v))| = o_p(1).$$

Further, note that

$$(40) \quad \sum_{\ell=1}^{9} Z_{0n}(\ell,u,v) = \frac{1}{n}\sum_{i=1}^{n}\mathbb{I}\{U_i \leq u, V_i \leq v\} + F_n^*(u) + G_n^*(v) + 1.$$

Combining (38), (39) and (40), implies that, uniformly in $(u,v)$,

$$\sqrt{n}(\hat{C}_n^{(\mathrm{MR})}(u,v) - C(u,v)) = B_n(u,v) + C_n^{h_n}(u,v) + o_p(1),$$

where $B_n$ is given by (25) and

$$C_n^{h_n}(u,v) = \sqrt{n}\left\{\sum_{\ell=1}^{9}\mathrm{E}\left[K\left(\frac{u-\hat{U}^{(\ell)}}{h_n}\right) - K\left(\frac{-\hat{U}^{(\ell)}}{h_n}\right)\right] \right.$$
$$\left. \times \left[K\left(\frac{v-\hat{V}^{(\ell)}}{h_n}\right) - K\left(\frac{-\hat{V}^{(\ell)}}{h_n}\right)\right] - C(u,v)\right\}.$$

STEP 2. Now, similarly as in Step 2 of the proof in Appendix A.1, we derive that, uniformly in $(u,v)$,

$$(41) \quad C_n^{h_n}(u,v) = -C_u(u,v)\mathbb{U}_n(u) - C_v(u,v)\mathbb{V}_n(v) + \sqrt{n}D_n(u,v) + o_P(1),$$

with the bias term $D_n(u,v)$ given by

$$(42) \quad D_n(u,v) = \sum_{\ell=1}^{9}\mathrm{E}\left[K\left(\frac{u - U_1^{(\ell)}}{h_n}\right) - K\left(\frac{-U_1^{(\ell)}}{h_n}\right)\right]$$
$$\times \left[K\left(\frac{v - V_1^{(\ell)}}{h_n}\right) - K\left(\frac{-V_1^{(\ell)}}{h_n}\right)\right] - C(u,v).$$

To show that this bias term is uniformly $O(h_n^2)$ is straightforward but tedious. The most simple case is if $(u,v) \in [h_n, 1-h_n]^2$. Then, (42) boils down to

$$D_n(u,v) = \mathrm{E} K\left(\frac{u-U_1}{h_n}\right) K\left(\frac{v-V_1}{h_n}\right) - C(u,v)$$

$$= \int_{-1}^{1} \int_{-1}^{1} C(u-th_n, v-sh_n) k(t) k(s) \, dt \, ds - C(u,v)$$

and the assertion follows simply by Taylor expansion.

Regarding the remaining cases, we will be dealing explicitly only with $(u,v) \in [1-h_n, 1]^2$. The other cases may be handled in a similar way.

Note that Taylor expansion together with the assumptions of the theorem imply $C(u,v) = u+v-1+O(h_n^2)$ uniformly in $(u,v) \in [1-2h_n, 1]^2$. Further, routine algebra shows that (42) simplifies to

$$\begin{aligned}
D_n(u,v) &= \mathrm{E} K\left(\frac{u-U_1}{h_n}\right) K\left(\frac{v-V_1}{h_n}\right) \\
&\quad + \mathrm{E} K\left(\frac{u+U_1-2}{h_n}\right) K\left(\frac{v-V_1}{h_n}\right) \\
&\quad + \mathrm{E} K\left(\frac{u-U_1}{h_n}\right) K\left(\frac{v+V_1-2}{h_n}\right) \\
&\quad + \mathrm{E} K\left(\frac{u+U_1-2}{h_n}\right) K\left(\frac{v+V_1-2}{h_n}\right) - C(u,v).
\end{aligned} \tag{43}$$

Let us compute

$$\begin{aligned}
&\mathrm{E} K\left(\frac{u-U_1}{h_n}\right) K\left(\frac{v-V_1}{h_n}\right) \\
&= \int_{(u-1)/h_n}^{1} \int_{(v-1)/h_n}^{1} C(u-th_n, v-sh_n) k(t) k(s) \, dt \, ds \\
&\quad + \int_{(u-1)/h_n}^{1} \int_{-1}^{(v-1)/h_n} (u-th_n) k(t) k(s) \, dt \, ds \\
&\quad + \int_{-1}^{(u-1)/h_n} \int_{(v-1)/h_n}^{1} (v-sh_n) k(t) k(s) \, dt \, ds \\
&\quad + \int_{-1}^{(u-1)/h_n} \int_{-1}^{(v-1)/h_n} 1 k(t) k(s) \, dt \, ds \\
&= \int_{(u-1)/h_n}^{1} \int_{(v-1)/h_n}^{1} (u-th_n+v-sh_n-1) k(t) k(s) \, dt \, ds \\
&\quad + \int_{(u-1)/h_n}^{1} (u-th_n) k(t) \, dt \, K\left(\frac{u-1}{h_n}\right)
\end{aligned} \tag{44}$$

$$+ \int_{(v-1)/h_n}^{1} (v - sh_n) k(s) \, ds \, K\left(\frac{v-1}{h_n}\right)$$

$$+ K\left(\frac{u-1}{h_n}\right) K\left(\frac{v-1}{h_n}\right) + O(h_n^2)$$

$$= \cdots$$

$$= (u + v - 1) + K\left(\frac{u-1}{h_n}\right)(1 - u) + K\left(\frac{v-1}{h_n}\right)(1 - v)$$

$$- h_n \int_{(u-1)/h_n}^{1} t k(t) \, dt - h_n \int_{(v-1)/h_n}^{1} t k(t) \, dt + O(h_n^2).$$

Similarly,

$$\mathrm{E} K\left(\frac{u + U_1 - 2}{h_n}\right) K\left(\frac{v - V_1}{h_n}\right)$$

(45)
$$= \int_{-1}^{(u-1)/h_n} \int_{(v-1)/h_n}^{1} P(U_1 > 2 + th_n - u,$$
$$V_1 \leq v - sh_n) k(t) k(s) \, dt \, ds$$

$$+ \int_{-1}^{(u-1)/h_n} \int_{-1}^{(v-1)/h_n} (1 - 2 - th_n + u) k(t) k(s) \, dt \, ds$$

$$= \cdots = (u - 1) K\left(\frac{u-1}{h_n}\right) - h_n \int_{-1}^{(u-1)/h_n} t k(t) \, dt + O(h_n^2),$$

(46)
$$\mathrm{E} K\left(\frac{v + V_1 - 2}{h_n}\right) K\left(\frac{u - U_1}{h_n}\right)$$
$$= (v - 1) K\left(\frac{v-1}{h_n}\right) - h_n \int_{-1}^{(v-1)/h_n} t k(t) \, dt + O(h_n^2),$$

(47)
$$\mathrm{E} K\left(\frac{u + U_1 - 2}{h_n}\right) K\left(\frac{v + V_1 - 2}{h_n}\right)$$
$$= O(h_n^2).$$

Combining (43) with (44), (45), (46) and (47) gives us

$$D_n(u, v) = u + v - 1 + O(h_n^2) - C(u, v) = O(h_n^2),$$

which was to be proved.

## APPENDIX B: PROOF OF THEOREM 2

**B.1. Weak convergence of the processes $\mathbb{C}_n^{(\mathrm{LLS})}$ and $\mathbb{C}_n^{(\mathrm{MRS})}$.** The proof of Theorem 2 for these estimators goes completely along the lines of the proof of Theorem 1, apart from a small difference in Step 2.

This difference is in calculating the bias term $D_n$ given by (35) for $\hat{C}_n^{(\text{LL})}$ and by (42) for $\hat{C}_n^{(\text{MR})}$. The shrinking of the bandwidths by the function $b(w) = \min\{\sqrt{w}, \sqrt{1-w}\}$, together with condition (9), guarantees that the Taylor expansion

$$C(u - sh_n b(u), v - th_n b(v))$$
$$= C(u,v) - sh_n b(u) C_u(u,v) - th_n b(v) C_v(u,v) + O(h_n^2)$$

holds uniformly in $(u,v) \in [0,1]^2$ and $(s,t) \in [-1,1]^2$. Applying the above expansion in the bias calculations completes the proof.

**B.2. Weak convergence of the process $\mathbb{C}_n^{(\mathbf{T})}$.** The proof is completely analogous to (and simpler than) the proof of Theorem 1 for $\hat{C}_n^{(\text{LL})}$. The only difference is in calculating the bias $D_n$, which is for the estimator $\hat{C}_n^{(\text{T})}$ given by

$$D_n(u,v) = \int_{-1}^{1} \int_{-1}^{1} C(\Phi(\Phi^{-1}(u) - sh_n), \Phi(\Phi^{-1}(v) - th_n)) k(s) k(t) \, ds \, dt.$$

As all the second order partial derivatives of $C(\Phi(\Phi^{-1}(u) - sh_n), \Phi(\Phi^{-1}(v) - th_n))$ taken as a function of $(s,t)$ are bounded by the assumptions of the theorem, Taylor expansion gives us $D_n(u,v) = O(h_n^2)$ uniformly in $(u,v)$ which proves the statement.

## APPENDIX C: JUSTIFICATION OF BOOTSTRAP TESTS BASED ON $\hat{C}_N^{(\text{LL})}$

Denote the process underlying the goodness-of-fit statistics as $G_n = \sqrt{n} \times (\hat{C}_n^{(\text{LL})} - C_{\hat{\theta}_n})$ and $G_n^* = \sqrt{n}(\hat{C}_n^{(\text{LL})*} - C_{\hat{\theta}_n^*})$ for its bootstrap version. In the following lemma, we will suppose that the true copula $C$ belongs to a known parametric family of copulas $\mathcal{C}_0 = \{C_\theta, \theta \in \Theta\}$.

LEMMA C.1. *Assume that the parametric family of copulas $\mathcal{C}_0$ satisfies the assumptions of Theorem 1 of Genest and Rémillard (2008) and, moreover, the derivative of the true copula $C_{\theta_0}$ with respect to $\theta$ is continuous as a function of $(u,v)$ in $[0,1]^2$. Then, there exists a "nonstochastic" sequence of functions $a_n$ from $[0,1]^2$ to $[0,1]$, such that $(G_n - a_n, G_n^* - a_n)$ converges in distribution to two independent copies of the same process.*

PROOF. From the proof of Theorem 1 of this paper, it follows that

$$\sqrt{n}\hat{C}_n^{(\text{LL})} = \sqrt{n}C_n^{(\text{E})} + \sqrt{n}D_n^{\theta_0} + o_p(1),$$

where

$$D_n^\theta(u,v) = \int_{-1}^{1} \int_{-1}^{1} [C_\theta(u - th_n, v - sh_n) - C_\theta(u,v)] k_{u,h}(s) k_{v,h}(t) \, dt \, ds.$$
(48)

Thus, the processes $G_n$ and $G_n^*$ might be rewritten as

(49) $$G_n = \sqrt{n}(C_n^{(\mathrm{E})} - C_{\hat{\theta}_n}) + \sqrt{n} D_n^{\theta_0} + o_p(1),$$

(50) $$G_n^* = \sqrt{n}(C_n^{(\mathrm{E})*} - C_{\hat{\theta}_n^*}) + \sqrt{n} D_n^{\hat{\theta}_n} + o_p(1).$$

As the empirical copula process $\sqrt{n}(C_n^{(\mathrm{E})} - C_{\theta_0})$ converges weakly, Theorem 1 of Genest and Rémillard (2008) implies that the first terms on the right-hand sides of (49) and (50) converge jointly in distribution to independent copies of the same process. Thus, defining $a_n(u,v)$ as $\sqrt{n} D_n^{\theta_0}(u,v)$, it remains to show that

$$\sup_{u,v} |D_n^{\theta_0}(u,v) - D_n^{\hat{\theta}_n}(u,v)| = o_p\left(\frac{1}{\sqrt{n}}\right).$$

But this follows directly from (48), the first order Taylor expansion of $C_{\hat{\theta}_n}$ around the true value of the parameter $\theta_0$ and the assumptions of the lemma. □

## APPENDIX D: VERIFICATION OF (9) FOR SOME FAMILIES OF COPULAS

We will verify assumption (9) only for $C_{uu}$. The assumptions about $C_{uv}$, $C_{vv}$ may be checked analogously.

**Clayton and Gumbel copulas.** Clayton and Gumbel copulas belong to an Archimedean family of copulas given by

(51) $$C(u,v) = \phi^{-1}(\phi(u) + \phi(v)),$$

where the function $\phi$ is called a generator of the copula. The generator of a Clayton copula is given by $\phi(t) = \frac{1}{\theta}(t^{-\theta} - 1)$ with $\theta \geq 0$ and that of a Gumbel copula by $\phi(t) = (-\log t)^\theta$ with $\theta \geq 1$.

Direct differentiaton of (51) yields

(52)
$$C_u(u,v) = \frac{\phi'(u)}{\phi'(C(u,v))},$$

$$C_{uu}(u,v) = \frac{\phi''(u)}{\phi'(C(u,v))} - \frac{[\phi'(u)]^2 \phi''(C(u,v))}{[\phi'(C(u,v))]^3}.$$

For a Clayton and a Gumbel copula, it is easy to verify that $\frac{\phi''(u)}{\phi'(u)} = O(\frac{1}{u(1-u)})$. Hence, we can bound the first term on the right-hand side of the expression for $C_{uu}(u,v)$ in (52) uniformly in $v$ by

$$\left|\frac{\phi''(u)}{\phi'(C(u,v))}\right| \leq \left|\frac{\phi''(u)}{\phi'(u)}\right|\left|\frac{\phi'(u)}{\phi'(C(u,v))}\right| \leq \left|\frac{\phi''(u)}{\phi'(u)}\right||C_u(u,v)|$$

(53)
$$= O\left(\frac{1}{u(1-u)}\right).$$

The second term on the right-hand side of the expression for $C_{uu}(u,v)$ in (52) is a more delicate one. For a Clayton copula, we have $\phi'(t) = -t^{-\theta-1}$ and $\phi''(t) = (\theta+1)t^{-\theta-2}$, which implies

$$\left|\frac{[\phi'(u)]^2\phi''(C(u,v))}{[\phi'(C(u,v))]^3}\right| = \left|\frac{(\theta+1)[C(u,v)]^{3\theta+3}}{u^{2\theta+2}[C(u,v)]^{\theta+2}}\right|$$

(54)
$$= \left|\frac{(\theta+1)[C(u,v)]^{2\theta+1}}{u^{2\theta+2}}\right|$$

$$\leq \frac{\theta+1}{u} = O\left(\frac{1}{u}\right),$$

using the Fréchet–Hoeffding upper bound for a copula [see Nelsen (2006)]. Combining (53) and (54) verifies (9) for $C_{uu}$ of a Clayton copula.

For a Gumbel copula we have

$$\phi'(u) = \theta(-\log u)^{\theta-1}\left(\frac{-1}{u}\right),$$

$$\phi''(u) = \theta(\theta-1)(-\log u)^{\theta-2}\left(\frac{1}{u^2}\right) + \theta(-\log u)^{\theta-1}\left(\frac{1}{u^2}\right),$$

which implies

$$\frac{[\phi'(u)]^2\phi''(C(u,v))}{[\phi'(C(u,v))]^3}$$

$$= \left(\theta^2(-\log u)^{2\theta-2}\left[\theta(\theta-1)(-\log C(u,v))^{\theta-2}\frac{1}{[C(u,v)]^2}\right.\right.$$

$$\left.\left.+ \theta(-\log C(u,v))^{\theta-1}\frac{1}{[C(u,v)]^2}\right]\right)$$

(55)
$$\times \left(u^2\left[\theta(-\log C(u,v))^{\theta-1}\frac{-1}{C(u,v)}\right]^3\right)^{-1}$$

$$= \frac{-(-\log u)^{2\theta-2}C(u,v)}{u^2}$$

$$\times [(\theta-1)(-\log C(u,v))^{1-2\theta} + (-\log C(u,v))^{2-2\theta}].$$

When $u \to 0_+$, the key fact is that

$$(56) \quad \frac{C(u,v)}{u^2} \frac{(-\log u)^{2\theta-2}}{(-\log C(u,v))^{2\theta-2}} \leq \frac{u}{u^2} \frac{(-\log u)^{2\theta-2}}{(-\log u)^{2\theta-2}} = \frac{1}{u},$$

and, when $u \to 1_-$,

$$(57) \quad \frac{(-\log u)^{2\theta-2}}{(-\log C(u,v))^{2\theta-1}} \leq \frac{(-\log u)^{2\theta-2}}{(-\log u)^{2\theta-1}} = \frac{1}{-\log u} = O\left(\frac{1}{1-u}\right).$$

Combining (53), (55), (56) and (56) verifies (9) for $C_{uu}$ of a Gumbel copula.

**Normal copula.** The normal copula is given by

$$C(u,v) = \int_{-\infty}^{\Phi^{-1}(u)} \int_{-\infty}^{\Phi^{-1}(v)} \frac{1}{2\pi\sqrt{1-\rho^2}} \exp\left\{\frac{s^2 - 2\rho st + t^2}{2(1-\rho^2)}\right\} ds\,dt,$$

$$\rho \in (-1,1),$$

where $\Phi$ is the cumulative distribution function of a standard normal variable.

By a direct computation (or with the help of properties of a conditional normal distribution), we get

$$C_u(u,v) = \Phi\left(\frac{\Phi^{-1}(v) - \rho\Phi^{-1}(u)}{\sqrt{1-\rho^2}}\right),$$

whose derivative with respect to $u$ is given by

$$(58) \quad C_{uu}(u,v) = \frac{-\rho}{\sqrt{1-\rho^2}} \phi\left(\frac{\Phi^{-1}(v) - \rho\Phi^{-1}(u)}{\sqrt{1-\rho^2}}\right) \frac{1}{\phi(\Phi^{-1}(u))},$$

where $\phi = \Phi'$. As $\phi$ is bounded, it is sufficient to deal with $[\phi(\Phi^{-1}(u))]^{-1}$. L'Hôpital's rule yields

$$\frac{u(1-u)}{\phi(\Phi^{-1}(u))} \sim \frac{1-2u}{\Phi^{-1}(u)} = o(1), \quad \text{for } u \to 0_+ \ (u \to 1_-),$$

which, together with (58), verifies (9) for $C_{uu}$ of a normal copula.

**Student copula.** The Student copula (with $m$ degrees of freedom) is given by

$$C(u,v) = \int_{-\infty}^{t_m^{-1}(u)} \int_{-\infty}^{t_m^{-1}(v)} \frac{1}{2\pi\sqrt{1-\rho^2}} \left(1 + \frac{s^2 - 2\rho st + t^2}{m(1-\rho^2)}\right)^{-(m+2)/2} ds\,dt,$$

$$\rho \in (-1,1),$$

where $t_m^{-1}(\cdot)$ is the quantile function of the Student distribution with $m$ degrees of freedom.

Direct calculation shows that

$$
\begin{aligned}
C_u(u,v) &= \frac{d^{(m+2)/2}}{f_m(t_m^{-1}(u))} \frac{1}{(c+d)^{(m+1)/2}} \\
&\quad \times \int_{-\infty}^{(t_m^{-1}(v)-\rho t_m^{-1}(u))/\sqrt{d+c}} (1+x^2)^{-(m+2)/2}\,dx,
\end{aligned}
\tag{59}
$$

where $c = [t_m^{-1}(u)]^2$, $d = m(1-\rho^2)$ and $f_m$ is the density of the Student distribution with $m$ degrees of freedom. Assumption (9) for $C_{uu}$ of a Student copula can be verified by differentiating (59) with respect to $u$. The useful facts (which follow by l'Hôpital's rule or properties of the density $f_m$) are

$$\frac{t_m^{-1}(u)}{u} \sim \frac{1}{f_m(t_m^{-1}(u))},$$

$$\frac{f_m'(t_m^{-1}(u))}{f_m(t_m^{-1}(u))} = O\left(\frac{1}{t_m^{-1}(u)}\right) \quad \text{for } u \to 0_+ \ (u \to 1_-),$$

where $f_m'$ is the derivative of $f_m$.

**Acknowledgments.** The authors are grateful to the co-Editor, an Associate-Editor and two reviewers for their very helpful comments, which led to a considerable improvement of the original version of the paper.

M. Omelka
Jaroslav Hájek Center for Theoretical
  and Applied Statistics
Charles University Prague
Sokolovská 83, 186 75 Praha 8
Czech Republic
E-mail: omelka@karlin.mff.cuni.cz

I. Gijbels
Department of Mathematics
  and Leuven Statistics
  Research Center (LStat)
Katholieke Universiteit Leuven
Celestijnenlaan 200B
B-3001 Leuven (Heverlee)
Belgium
E-mail: irene.gijbels@wis.kuleuven.be

N. Veraverbeke
Center for Statistics
Hasselt University
Agoralaan—building D
B-3590 Diepenbeek
Belgium
E-mail: noel.veraverbeke@uhasselt.be